\documentclass[11pt,           
english                       
]{article}

\usepackage{stmaryrd}
\usepackage[english,german,strings]{babel}     
\usepackage{amsmath}           
\usepackage[utf8]{inputenc}    
\usepackage[T1]{fontenc}       
\usepackage{longtable}         
\usepackage{exscale}           
\usepackage[final]{graphicx}   
\usepackage{array}             
\usepackage{fancyhdr}          
\usepackage[a4paper]{geometry} 
\usepackage[multiuser]{fixme}  
\usepackage{xspace}            
\usepackage{tikz}              
\usepackage{appendix}
\usepackage[expansion=false    
]{microtype}        
\usepackage[nottoc]{tocbibind}
\usepackage[backref=page,      
final=true,         
pdfpagelabels,       
hidelinks 
]{hyperref}         
\usepackage[shortlabels]{enumitem}

\usepackage[numbers]{natbib}
\usepackage{tikz-cd}
\usepackage{tikz} 
\usepackage{nchairx} 
\usepackage{soul}
\usepackage{combelow}

\newcommand{\blup}{\mathrm{Blup}}
\newcommand{\Partial}[1]{\frac{\partial}{\partial #1}}

\newcommand{\blowdown}[3]{#1 \colon \blup(#2,#3) \to #2 }
\newcommand{\foliation}[1]{\mathcal{#1}}

\newcommand{\vanishing}{\mathcal{I}}

\newcommand{\lin}{\mathrm{lin}}
\newcommand{\sdiv}{\mathbb{P}(\nu_N(M))}

\newcommand{\RR}{\ensuremath{\mathbb R}}
\newcommand{\TT}{\ensuremath{\mathbb T}}
\newcommand{\PP}{\ensuremath{\mathbb P}}

\newcommand{\g}{\ensuremath{\mathfrak{g}}}
\newcommand{\h}{\ensuremath{\mathfrak{h}}}

\newcommand{\cO}{\mathcal{O}} 
\newcommand{\diffto}{\xrightarrow{\raisebox{-0.2 em}[0pt][0pt]{\smash{\ensuremath{\sim}}}}}

\usetikzlibrary{cd}
\usepackage{amsmath, amssymb}

\newcommand{\plin}{\pi_{\lin}}

\newcommand{\DEC}{\mathrm{d}_{\mathfrak{g}}}

\newcommand{\DECC}{\mathrm{d}_{\mathfrak{g}_{\mathbb{C}}}}

\usepackage{color} \definecolor{tocolor}{rgb}{.1,.1,.5}
\definecolor{urlcolor}{rgb}{.2,.2,.6}
\definecolor{linkcolor}{rgb}{.1,.1,.6}
\definecolor{citecolor}{rgb}{.6,.2,.1}
\hypersetup{colorlinks=true, urlcolor=urlcolor,
  linkcolor=linkcolor, citecolor=citecolor}
\definecolor{darkgreen}{rgb}{0.0, 0.5, 0.0}

\begingroup
\hypersetup{linkcolor=blue}
\endgroup

\begin{document}
\selectlanguage{english}
\title{Blowups of Dirac structures}

\author{
Ioan M\u{a}rcu\cb{t} 
\thanks{\textsf{Mathematisches Institut, Universit\"at zu K\"oln,
Weyertal 86-90
50931 K\"oln, Germany. 
}}
\and
Andreas Sch\"{u}\ss ler
\thanks{\textsf{Department of Mathematics, KU Leuven, Celestijnenlaan 200B, 3001 Leuven, Belgium. \\
\indent Email: \texttt{andreas.schuessler@kuleuven.be}}}
  \and
 Marco Zambon
\thanks{\textsf{Department of Mathematics, KU Leuven, Celestijnenlaan 200B, 3001 Leuven, Belgium.\\
\indent Email: \texttt{marco.zambon@kuleuven.be}}}
}

\date{}
\maketitle
  


\begin{abstract}
Given a real, twisted Dirac structure $L$ on a smooth manifold $M$, and a closed embedded submanifold $N\subseteq M$ of codimension $>1$, we characterise when $L$ lifts to a smooth, twisted Dirac structure on the real projective blowup of $M$ along $N$. This holds precisely when $N$ is either a submanifold transverse to $L$ (with no further restrictions) or a submanifold invariant for $L$, for which the Lie algebras transverse to $N$ have all of the same constant height $k\geq 0$. In this paper, we also classify Lie algebras satisfying this Lie-theoretic property up to isomorphism: Lie algebras of constant height $k=0$ are either abelian or a semi-direct product $\mathbb{R}\ltimes \mathbb{R}^n$ for the diagonal representation of $\mathbb{R}$ on $\mathbb{R}^n$; there is only $\mathfrak{so}(3)$ of constant height $k=1$; there are no Lie algebras of constant height $k\geq 2$. We recover a theorem of Polishchuk, which establishes that a Poisson structure lifts to a Poisson structure on the blowup of a submanifold exactly when the submanifold is invariant and the transverse Lie algebras have constant height $k=0$. 
 \end{abstract}
\setcounter{tocdepth}{1} 
\tableofcontents

\section{Introduction}

 Given a manifold $M$ and a closed and embedded submanifold $N$, the \emph{real projective blowup} is obtained by replacing $ N $ by $\sdiv$, the projectivisation of the normal bundle $\nu_N(M)\to N$ of $ N $ in $M$. This yields a smooth manifold (without boundary), which we denote by $\blup(M,N)$. The blowup comes with a canonical blowdown map  $$p\colon \blup(M,N)\to M,$$  which  
on $\field P(\nu_N(M))=p^{-1}(N)$ is the natural projection to $N$, and which defines a diffeomorphism from the complement of $\field P(\nu_N(M))$ to the complement of $N$. In particular, $p$ is not a submersion.

One can address the question of when geometric structures on $M$ ``lift'' to the blowup.
This question has been addressed in many contexts in the literature, both to produce new examples and as a way to desingularise geometric structures.

Notice that we are considering real blowups,
and that in the literature complex blowups also appear frequently. 

In \cite{polishchuk:1997a}, Polishchuk studied the projective blowup of Poisson schemes, i.e.\ the existence of a Poisson structure on the blowup such that the blowdown map is Poisson. 
We state the result of \cite[\S 8]{polishchuk:1997a} in its weaker, differential-geometric formulation, in a similar way to \cite[\S 2.2]{gualtieri.li:2014a}.

Recall that, given a Poisson submanifold $N$ of a Poisson manifold $(M,\pi)$, the conormal space at any $q\in N$,
\begin{equation*}
    (T_qN)^\circ := \ann (T_qN)=(T_qM/T_qN)^\ast
\end{equation*}
carries a Lie bracket, turning $(TN)^\circ$ into a bundle of Lie algebras. Explicitly, the bracket is given by \begin{equation}\label{eq:Liealgconormal}
    [(\D f)(q), (\D g)(q)]:= (\D \{f,g\})(q),
 \end{equation}
where $f,g$ are smooth functions on $M$ vanishing on $N$.

\begin{nntheorem}[Polishchuk \cite{polishchuk:1997a}]
    Let $(M,\pi)$ be a Poisson manifold and $N\subseteq M$ a closed and embedded Poisson submanifold. There exists a Poisson bivector field $\tilde{\pi}$ on $\blup(M,N)$ which is $p$-related to $\pi$ 
    if and only if each Lie algebra $
    (T_qN)^{\circ}$, for $q\in N$, is either abelian or isomorphic to the semi-direct product by the diagonal representation $\lambda\mapsto \lambda \id_{\field R^n}$ of $\field R$ on $\field R^n$, denoted by $\field R\ltimes \field R^n$.
\end{nntheorem}

In this paper we shall be interested in lifting \emph{twisted Dirac structures}. Dirac structures include  Poisson structures, closed 2-forms, and foliations as special cases, and are suitable to describe geometric structures that become ``infinite'' at certain points.

Given a (twisted) Dirac structure $L\subseteq TM\oplus T^\ast M$ on a manifold $M$ and a closed, embedded submanifold $N\subseteq M$ with $\codim(N)> 1$, we ask:\\

\emph{Does there exist a (twisted) Dirac structure $\tilde{L}$ on $\blup(M,N)$, such that}
\begin{equation*}
\tilde{L}\at{\blup(M,N)\setminus \field P(\nu_N(M))}=p^* (L\at{M\setminus N})\text{?}
\end{equation*}
If such a Dirac structure $\tilde{L}$ exists, it is necessarily unique. In this case, we say that $L$ \emph{lifts to the blowup}. Further, we ask whether \emph{the blowdown map is forward or backward Dirac}.

Our main theorem is an extension of Polishchuk's result (in the smooth setting, formulated before) to twisted Dirac structures. 
To state it, first notice that, given a twisted Dirac structure $L$ on $M$ and a submanifold $N\subseteq M$ which is invariant (i.e.\ a union of leaves), the annihilator $(TN)^\circ$ is a bundle of Lie algebras, as in the case of Poisson structures. Indeed, $(TN)^\circ$ fits into a short exact sequence of Lie algebroids:
\[0\longrightarrow (TN)^{\circ}\longrightarrow L\at{N}\longrightarrow \iota_N^!L\longrightarrow 0,
\]
where $L\at{N}$ denotes the pullback of the Lie algebroid $L$ to $N$ and $\iota_N^!L$ denotes the pullback of the twisted Dirac structure $L$ to $N$.

\begin{nnmaintheorem} Let $ L $ be a twisted Dirac structure on a manifold $M$ and $N\subseteq M$ a connected, closed, and embedded submanifold of codimension $\codim N>1$. Then, $L$ lifts to a  twisted Dirac structure on $ \blup(M,N) $ if and only if one of the following conditions holds.
	\begin{enumerate}[(1)]
    \item $N\subseteq M$ is a transversal.
    \item $N\subseteq M$ is an invariant submanifold and the bundle of Lie algebras $(TN)^{\circ}\to N$ satisfies one of the following conditions:
    \begin{enumerate}[(a)]
        \item Each fibre is either abelian or isomorphic to the semi-direct product $ \field R\ltimes \field R^n$ by the diagonal representation of $\field R$ on $\field R^n$,
        \item The fibres are  all isomorphic to  $\liealg{so}(3)$. 
    \end{enumerate}
    \end{enumerate}
 The blowdown map is a backward Dirac map in case (1), and it is a forward Dirac map in case (2).

\end{nnmaintheorem}
Invariant and transversal submanifolds are defined in Definition \ref{def:invtrans}. Note that when $\codim(N)=1$ one has  $ \blup(M,N)\simeq M$ via the blowdown map, so and Dirac structure lifts to $\blup(M,N)$.

\begin{nnremark}
If $N$ is a presymplectic leaf of the Dirac structure $L$, or equivalently, if $N$ is invariant and $L\at{N}$ is a transitive Lie algebroid, then $(TN)^{\circ}$ is a locally trivial bundle of Lie algebras. Hence all fibres are isomorphic to each other, simplifying the statement of the main theorem.
\end{nnremark}
\begin{nnexample}
    As a simple illustration of the above theorem, let $\g$ be a Lie algebra and $\h\subseteq \liealg g$ a proper ideal. Then $\h^{\circ}$ is a Poisson submanifold of $\g^*$, endowed with the  standard linear Poisson structure. The latter lifts to a Dirac structure on $\blup(\g^*,\h^{\circ})$ if and only if the Lie algebra $\h$ is isomorphic to one of the Lie algebras appearing in item (2) of the main theorem.
\end{nnexample}

\noindent{\bf Outline of the paper and intermediate results.} Sections \ref{sec:blowup} and \ref{sec:Dirac} are devoted to recalling basic notions and results about the blowup construction and Dirac structures, respectively.

The proof of the Main Theorem, is divided into several steps corresponding to the remaining sections of the paper. 

The case of transverse submanifolds is easily dealt with in 
Section \ref{sec:blowup_trans}, where the following is proven.

\begin{nntheorem}[\ref{thm:lift_to_transversals}] 
Let $L$ be a twisted Dirac structure on $M$. For any closed, embedded trans\-verse submanifold $N\subseteq M$, $L$ lifts to $\blup(M,N)$ and the blowdown map is a backward Dirac map.  
\end{nntheorem}

In Section \ref{section:dichotomy} we prove the transverse or invariant dichotomy, from the Main Theorem.

\begin{nntheorem}[\ref{theorem:Dirac_everything}] 
  Let $L$ be a twisted Dirac structure on $M$ and $N\subseteq M$ a closed, embedded, and connected submanifold of codimension $>1$. If $L$ lifts to a twisted Dirac structure on $\blup(M,N)$, then one of the two conditions holds.
    \begin{itemize}
        \item $N\subseteq M$ is a transversal. In this case, the blowdown map is a backward Dirac map.
        \item $N\subseteq M$ is an invariant submanifold.
        In this case, the blowdown map is a forward Dirac map.
    \end{itemize}
\end{nntheorem}

To prove this result, we make use of the description of Dirac structures in terms of spinor lines, i.e.\ line subbundles of $\oplus_k \wedge^k T^*M$ \cite{gualtieri:2011a}. For this we develop in Subsection \ref{section:extending_line_bundles} techniques to study the problem of extending line bundles along codimension 
 one submanifolds---some of the  results obtained are similar those in \cite{blohmann:2017a}. In particular, our techniques imply the following result. 
 
\begin{nncorollary}[\ref{coro:almost_everywhere}]
   Let $L$ be a twisted Dirac structure on $M$ and $N\subseteq M$ a closed and embedded submanifold. There exists an open and dense subset
   $V\subseteq \field P(\nu_N(M))$ such that the Dirac structure $p^\ast(L\at{M\setminus N})$
    extends smoothly to a Dirac structure on $\big(\blup(M,N)\setminus \field P(\nu_N(M))\big) \cup V$.
\end{nncorollary}

In Section \ref{section:Poisson_structures_by_spinors} we turn to the general problem of blowing up invariant submanifolds. 

The tools developed in Section \ref{section:dichotomy} are crucial also here. Namely, we use spinors to describe Dirac structures, and study the vanishing order along the divisor of the pullback of the spinor to the blowup.

Our goal is to reduce the problem to the following 
purely Lie theoretical property of the isotropy Lie algebras $(T_qN)^{\circ}$, $q\in N$. For a Lie algebra $(\g,[\argument,\argument])$, define the \emph{height} of an element $\xi \in \g^*\setminus\{0\}$ as 
  the integer 
$k\in \field N_0$ satisfying 
\begin{equation*}
        \xi\wedge (\DEC\xi)^{k}\neq 0 \quad\text{ and }\quad \xi\wedge (\DEC\xi)^{k+1}= 0,
\end{equation*}
where $\DEC$ denotes the Chevalley-Eilenberg differential. We say that $\mathfrak{g}$ is a Lie algebra of \emph{constant height} $k\geq 0$, if all elements $\xi \in \g^*\setminus\{0\}$ have height $k$.

Next, we note that the structure of a bundle of Lie algebras of $(TN)^{\circ}$ induces a fibrewise linear Poisson structure on the normal bundle $\nu_N(M)\simeq ((TN^{\circ}))^{*}$, which we denote by $\plin$. The Dirac structure $\mathrm{graph}(\plin)$ can be thought of as the linearization of the vertical component of $L$ at $N$.

The main results of Section \ref{section:Poisson_structures_by_spinors} are summarised in the following.

\begin{nntheorem}[\ref{theorem:lifting_by_spinors}, \ref{coro:vanishing:order}, \ref{pullback of invariant}, \ref{theorem:summary_submanifold}] Let $L$ be a twisted Dirac structure on $M$ and $N\subseteq M$ a closed, embedded and connected submanifold, which is invariant. The following are equivalent.
\begin{itemize}
\item $L$ lifts to  $\blup(M,N)$.
\item There exists $\ell\in \{1,\ldots, \mathrm{codim}\, N-1\}$ such that, for any (local) spinor $\phi$ for $L$, $p^*\phi$ has constant vanishing order $\ell$ along $\sdiv$.
\item  $\graph(\plin)$ lifts to $\blup(\nu_N(M),0_N)$.
\item The Lie algebras $(T_qN)^{\circ}$, $q\in N$, all have the same constant height $k$.
\end{itemize}
\end{nntheorem}

Finally, using the structure theory of semisimple Lie algebras, in Section \ref{sec:classif} we classify Lie algebras of constant height.

\begin{nntheorem}[\ref{theorem:classification_LA_of_constant_height}]
  Any Lie algebra $\liealg g$ of constant height is isomorphic to one of the following.
\begin{itemize}        \item An abelian Lie algebra $\field R^n$---this has height $0$.
 \item The semi-direct product $\field R\ltimes \field R^n$, for the representation $\lambda\mapsto \lambda\id_{\RR^n}$---this has height $0$.
        \item The Lie algebra $\liealg{so}(3)$---this has height $1$.
    \end{itemize}
\end{nntheorem}

Theorems 
\ref{thm:lift_to_transversals}, \ref{theorem:Dirac_everything}, \ref{theorem:lifting_by_spinors} and \ref{theorem:classification_LA_of_constant_height} together yield the Main Theorem.

In Section \ref{ref:sec:lift_Poisson},
for the special case of blowing up a zero of a Poisson structure, we provide an alternative, more geometric proof of Theorem \ref{theorem:lifting_by_spinors}, without using spinors, independently of our previous arguments. We do so in
Theorem \ref{theorem:point_lift}, which gives several geometric characterisations of the existence of lifts, one of them in terms of the dimensions of the coadjoint orbits.

\bigskip
\noindent{\bf Relation to the literature.} 
The main question we address in this paper is motivated by Polishchuk's result \cite[\S 8]{polishchuk:1997a} recalled above. In particular, our Main Theorem extends Polishchuk's result in the framework of smooth manifolds from Poisson to Dirac structures. Our techniques differ from those used by Polishchuk in \cite[\S 8]{polishchuk:1997a}. There, Poisson structures are described by means of the Poisson bracket on the algebra of functions; since for Dirac structures only the subalgebra of admissible functions is endowed with a Poisson bracket, we use more geometric techniques. It is however surprising that, although we allow for completely general submanifolds and completely general twisted Dirac structures, we obtain only a few new situations compared with Polishchuk's setting in which the Dirac structure lifts to the blowup: transverse submanifolds and invariant submanifolds of codimension three with transverse Poisson structure $\liealg{so}(3)$.

Indirectly, the case of $\liealg{so}(3)$ has appeared in the literature in relation to the problem of linearising a Poisson structure $\pi$ around a zero $q$ (i.e.\ $\pi(q)=0$). Namely, if the isotropy Lie algebra of $\pi$ at $q$ is $\liealg{so}(3)$, a geometric proof of Conn's linearisation theorem \cite{conn85} can be obtained by first blowing up $q$, then using Reeb's stability for the regular foliation underlying the lifted Dirac structure, and finally using a foliated Moser trick. This proof, envisioned by Weinstein \cite{weinstein83a}, was carried out by Dazord in \cite{dazord83}, and in fact, using the language of completely integrable 1-forms, the argument can be traced back to Reeb's thesis \cite{reeb52}.

Blowups have also been studied in the setting of generalised complex geometry in \cite{cavalcanti.gualtieri09} in dimension 4 and in \cite{bailey.cavalcanti.vdleerduran:2019a} in greater generality. The natural operation in this setting is the complex blowup, arising from a complex structure on the normal bundle to the submanifold. Like a foreshadowing of the ``invariant-transverse-dichotomy'' revealed in our paper, the submanifolds considered in \cite{bailey.cavalcanti.vdleerduran:2019a} are certain classes of Poisson submanifolds for which Polishchuk's result applies and certain Poisson transversals, with a transverse complex structure. In the second case, in order to obtain a generalised complex structure on the blowup, the lifted (complex) Dirac structure needs to be modified around the divisor using an adaptation of blowup construction in symplectic geometry \cite{guillemin.sternberg89a} and Lerman's symplectic cut construction \cite{lerman95a} to generalized complex geometry.

It would be interesting to see if there is more flexibility in the more general setting of weighted blowups, which have recently proven to be very successful in simplifying Hironaka's resolutions of singularities in algebraic geometry \cite{abramovich.temkin.wlodarczyk,mcguillan}. For Poisson structures, this is currently being investigated by Pym and collaborators \cite{WeightedPoisson}. 
 
Every Dirac manifold $(M,L)$ carries a natural singular foliation $\mathcal{F}=\mathrm{pr}_{TM}(\Gamma_c^\infty(L))$, whose leaves are precisely the presymplectic leaves. For an invariant (respectively transverse) submanifold $N\subseteq M$,  the singular foliation $\mathcal{F}$ lifts to a singular foliation $p^*\mathcal{F}$ on the blowup $\blup(M,N)$---see \cite[\S 1.5.7]{singfolnotes} (respectively 
\cite[\S 1.5.4]{singfolnotes}). In the case when the Dirac structure also lifts, one has the inclusion 
$p^*\mathcal{F}\subseteq \widehat{\mathcal{F}}$ 
into the singular foliation $\widehat{\mathcal{F}}$ of the lifted Dirac structure, which in general is strict.
When $N$ is a transverse submanifold, the above inclusion is an equality \cite{schuessler.zambon:2025a}.

Crainic, Fernandes, and  Mart\'inez Torres \cite{PMCT3} develop a procedure to lift a Poisson structure of compact type on a manifold $M$ to a regular Dirac structure on another manifold $\widehat{M}$ of the same dimension. However, $\widehat{M}$ is in general not a real projective blowup. In the special case of the Poisson manifold $\liealg{so}(3)^*$, their construction does agree with the real projective blowup, and fits in the setting of our Main Theorem and of Theorem \ref{theorem:point_lift}.
 
Finally, Theorem \ref{theorem:lifting_by_spinors} classifies (real) Lie algebras of constant height. We believe that this statement and its proof are of independent interest. In Remark \ref{rem:gozeremm} we relate the height of an element with the Cartan class of that element, as used for instance in \cite{goze.remm:2019a} to obtain classification results.

\bigskip
\noindent{\textbf{Acknowledgements.}} 
We thank Sergei Berezin for interesting discussions.

M.Z.\,acknowledges partial support by
the FWO and FNRS under EOS project G0I2222N, by FWO project G0B3523N, and by Methusalem grant METH/21/03 - long term structural funding of the Flemish Government (Belgium).

A.S.\,was financially supported by Methusalem grant METH/21/03 – long term structural funding of the Flemish Government (Belgium).

\section{Real projective blowups}\label{sec:blowup}

In this section, we briefly review the real projective blowup and some of its basic properties. 

As a set, the blowup of a closed, embedded submanifold $N$ of a manifold $M$ is given by
\begin{equation*}
	\blup(M,N)=(M\setminus N)\sqcup \field P(\nu_N(M)),
	\end{equation*}
	i.e.\,$ N $ is replaced with the projectivisation of the normal bundle $\nu_N(M):=TM\at{N}/TN$ of $ N $ in $ M $. 
    The blowup comes with a canonical smooth structure, for which the blowdown map,
    \begin{equation*}
        p\colon \blup(M,N)\to M,
    \end{equation*}
    is smooth. This map is the identity on $M\setminus N$ and the projection on $\field P(\nu_N(M))\to N$. If $ \codim N=1 $, then $p$ is a diffeomorphism. The smooth structure can be described via the following set of charts. 

   \subsubsection*{Charts for $ \blup(M,N) $}
		Let $ (U,(x,y)) $ be a submanifold chart for $N$, i.e.\ $ U\cap N=\{ x=0 \} $. Then the collection \begin{equation*}
		U_i= p^{-1}\big(\{ x_i\neq 0 \}\big)\cup \big\{ [v]\in \field P(\nu_{N}(M))\at{U\cap N}) \colon \D x_i(v)\neq 0  \big\},\qquad i=1,\dots,\codim(N),
		\end{equation*}
		is an open cover of $ p^{-1}(U) $. On each $U_i$, one defines coordinates $(\tilde{x},y)$ on $\blup(M,N)$, in which the blowdown map reads (see e.g.\ \cite[Remark 5.29]{obster:2021a} for details)
\begin{equation}\label{eq:blowdowncoord}
p(\tilde{x}_1,\dots,\tilde{x}_i,\dots,\tilde{x}_{\codim(N)},y)=(\tilde{x}_i\tilde{x}_1,\dots,\tilde{x}_i,\dots,\tilde{x}_i\tilde{x}_{\codim(N)},y).
        \end{equation}
Notice that the hyperplane $\field P(\nu_N(M))\cap U_i$ is given by $\tilde{x}_i=0$, and that the chart obtained by restriction is the well-known chart on the projective bundle induced by the fibrewise linear coordinates $(\D x,y)$ on the normal bundle $\nu_N(M)$.

\subsubsection*{Lifts of vector fields to $ \blup(M,N) $}

In the chart $(U_i, (\tilde{x},y))$, the lifts of the coordinate vector fields $\Partial{x_i}$ and $\Partial{x_k}$ for $k\neq i$ are given by 
\begin{equation}\label{eq:vfliftcoords}
\Partial{\tilde{x}_i}- \sum_{k\neq i} \frac{ \tilde{x}_k}{\tilde{x}_i}\Partial{\tilde{x}_k}\quad\text{ and }\quad \frac{ 1}{\tilde{x}_i}\Partial{\tilde{x}_k},
		\end{equation}
respectively, as one sees by applying these vector fields to functions of the form $p^\ast x_j$. This implies the following standard result (see e.g.\ \cite[Prop.\ 1.5.40]{singfolnotes} for points or \cite[Lemma 3.5]{schuessler:2024a}).
\begin{lemma}\label{lemma:lifting_vf_to_the_blowup}
		A vector field $X$ on $M$ is $p$-related to some vector field $\tilde{X}$ on $\blup(M,N) $ if and only if $ X $ is tangent to $ N $. In that case, $ \tilde{X} $ is unique and tangent to $ \field P(\nu_{N}(M))\subseteq \blup(M,N) $. 
	\end{lemma}

\section{Dirac structures}\label{sec:Dirac}

We give a short introduction to Dirac structures by first describing them as subbundles of the standard Courant algebroid and then via their spinor lines. Dirac structures were introduced in \cite{courant:1990a, courant.weinstein:1988a} and we recommend \cite{bursztyn:2013a} for a brief introduction to the field, and \cite{gualtieri:2011a} for the use of spinors.

\subsection{Dirac structures as Lagrangian subbundles of $\mathbb{T}M$}

Let $M$ be a manifold and $H\in \Omega^3(M)$ a closed $3$-form. The \emph{generalised tangent bundle} 
\begin{equation*}
    \mathbb{T}M:=TM\oplus T^\ast M
\end{equation*}
carries a \emph{nondegenerate pairing} and the \emph{Dorfman bracket}, defined by
\begin{equation*}
    \begin{aligned}
        \langle \argument,\argument \rangle\colon \mathbb{T}M\times \mathbb{T}M&\to \field R\\
        \begin{pmatrix}v\\
        \xi   \end{pmatrix},   \begin{pmatrix}w\\
        \eta   \end{pmatrix}&\mapsto \frac{1}{2}(\xi(w)+\eta(v))
    \end{aligned}
\end{equation*}
\begin{equation*}
    \begin{aligned}        [\argument,\argument]\colon\Gamma^\infty(\mathbb{T}M)\times \Gamma^\infty(\mathbb{T}M)&\to \Gamma^\infty(\mathbb{T}M)\\
        \begin{pmatrix}X\\
        \alpha   \end{pmatrix},\begin{pmatrix}Y\\
        \beta   \end{pmatrix}&\mapsto \begin{pmatrix}[X,Y]\\
        \Lie_X \beta-\mathrm{i}_Y \D \alpha+ \I_Y \I_X H\end{pmatrix}.
    \end{aligned}
\end{equation*}
Together with the \emph{anchor} map $\sharp=\pr_{TM}\colon \mathbb{T}M\to TM$ this endows $\mathbb{T}M$ with the structure of a \emph{$H$-twisted Courant algebroid}. 

\begin{definition}
    A subbundle $L\subseteq \mathbb{T}M$ is called an \emph{$H$-twisted Dirac structure} on $M$ if
    \begin{enumerate}
        \item $L$ is maximally isotropic (i.e.\ \emph{Lagrangian}) with respect to the pairing $\langle\argument,\argument\rangle$,
        \item $L$ is involutive, i.e.\ $[\Gamma^\infty(L),\Gamma^\infty(L)]\subseteq \Gamma^\infty(L)$.
    \end{enumerate}
When we do not  want to specify $H$, we call $L$ a simply a \emph{twisted Dirac structure}, or even simply \emph{a Dirac structure} (even if $H\neq 0$).
\end{definition}

Here are the standard examples for $H$-twisted Dirac structures, following \cite[Example 2.11-2.13]{gualtieri:2011a}.
\begin{enumerate}
    \item The graph
    of a bivector field $\pi$ on $M$,
    \begin{equation*}
        \graph(\pi)=
        \left\{\begin{pmatrix}\pi^\sharp \xi\\
        \xi\end{pmatrix} \in \mathbb{T}M\colon  \xi\in T^\ast M  \right\},
    \end{equation*} if and only if 
    \begin{equation*}
\Schouten{\pi,\pi}=2(\wedge^3\pi^{\sharp})H, 
    \end{equation*}

where $\Schouten{\argument,\argument }$ is the Schouten bracket. Then $\pi$ is called an \emph{$H$-twisted Poisson structure} \cite{severa.weinstein:2001a}.
    \item The graph of a $2$-form $\omega$ on $M$,
    \begin{equation*}
        \graph(\omega)=\left\{\begin{pmatrix} v\\
        \omega^\flat v\end{pmatrix} \in \mathbb{T}M\colon v\in T M\right\},
    \end{equation*}
    if and only if $\D\omega=H$.
    \item The subbundle $T\foliation{F}\oplus (T\foliation{F})^\circ$, for a regular foliation $\foliation F$ on $M$, if and only if $H\at{T\foliation{F}}=0$. Here, $\argument^\circ$ denotes the annihilator.
\end{enumerate}

\begin{definition}[\cite{bursztyn.radko:2003a, weinstein:1982a}]
Let $f\colon N\to M$ be a smooth map. Let $L_N$ be an $H_N$-twisted Dirac structure on $N$, and $L_M$ be an $H_M$-twisted Dirac structure on $M$. If $H_N=f^*H_M$, we call $f$
\begin{enumerate}
    \item \emph{a forward Dirac map}, if for all $n\in N$ we have $(L_M)_{f(n)}=f_!((L_N)_{n})$, where 
	\begin{equation*}
		f_!((L_N)_{n}):=\left\{ \begin{pmatrix}(T_{n}f) w\\
            \xi\end{pmatrix}
            \in \field T_{f(n)}M\colon 
            \begin{pmatrix} w \\(T_nf)^\ast\xi \end{pmatrix}\in (L_N)_n\right\};
	\end{equation*}
    \item \emph{a backward Dirac map}, if for all $n\in N$ we have $(L_N)_{n}=f^!((L_M)_{f(n)})$, where 
    \begin{equation*}
        f^!((L_M)_{f(n)})=\left\{  \begin{pmatrix} w \\(T_nf)^\ast\xi \end{pmatrix}
            \in \field T_{n}N\colon 
            \begin{pmatrix} (T_n f)w \\\xi \end{pmatrix}\in (L_M)_{f(n)}\right\}.
    \end{equation*}
\end{enumerate}
\end{definition}
Forward Dirac maps generalise Poisson maps, while backward maps generalise pullbacks of $2$-forms.

The \emph{gauge action} of a 2-form $B\in \Omega^2(M)$ on a subbundle $L\subseteq \field T M$ is defined by \cite{severa.weinstein:2001a}
    \begin{equation*}
        \E^B L := \left\{ \begin{pmatrix}X\\
            \alpha+\I_X B\end{pmatrix}
            \in \field TM \colon 
            \begin{pmatrix}X\\\alpha \end{pmatrix}\in L\right\}.
    \end{equation*}
Gauge actions send Lagrangian subbundles to Lagrangian subbundles. Moreover, $L\subseteq \field T M$ is an $H$-twisted Dirac structure if and only if $\E^B L \subseteq \field T M$ is an $(H+ \D B)$-twisted Dirac structure. 

\begin{definition}\label{def:isomorphic_dirac_structures} For $i=1,2$, let $L_i$ be an $H_i$-twisted Dirac structure on $M_i$. We say that $L_1$ and $L_2$ are \emph{isomorphic} if there exists a diffeomorphism $f\colon M_1\to M_2$ and a 2-form $B_2\in \Omega^2(M_2)$ such that 
    \begin{equation*}
        L_1=f^\ast(\E^{B_2}L_2) \quad\text{ and }\quad H_1=f^\ast(H_2+\D B_2).
    \end{equation*}
    We call $L_1$ and $L_2$ \emph{locally isomorphic} around $x_1\in M_1$ and $x_2\in M_2$ if their restrictions to open neighbourhoods $U_1$ of $x_1$ and $U_2$ of $x_2$, respectively, are isomorphic.
\end{definition}

\subsection{Spinor description of Dirac structures}\label{sec:spinor}

Another way to describe Dirac structures is by means of pure spinor lines.
Let $ M$ be a manifold. Define an action $\rho$ of $ \Gamma^\infty(\field TM) $ on $ \Omega^\bullet(M) $ by
\begin{equation*}	\rho\begin{pmatrix}X\\\alpha\end{pmatrix}\phi:= \I_X\phi+\alpha\wedge \phi.
\end{equation*}
Then, for a nowhere vanishing $ \phi\in \Omega^\bullet(M) $, called a \emph{spinor}, the spaces 
\begin{equation*}
	L_\phi:=\left\{ \begin{pmatrix}X\\\alpha\end{pmatrix}\in \field TM\colon \rho\begin{pmatrix}X\\\alpha\end{pmatrix}\phi=0 \right\}
\end{equation*}
are isotropic in every point of $ M $, i.e.\ $ L_\phi\subseteq L_\phi^\perp $, and $ \phi $ is called a \emph{pure spinor} if these spaces are of maximal dimension, i.e.\ $\rank L_\phi=\dim M$. 
The maximal isotropic subbundle corresponding to a pure spinor depends only on the line bundle 
\[\Sigma:=\mathbb{R}\cdot \phi\subseteq \wedge^{\bullet} T^*M.\] 
Conversely, every maximal isotropic subbundle has a corresponding pure spinor line bundle (\cite[III.1.9]{chevalley:1997a}, \cite[Proposition 1.3]{gualtieri:2011a}). 
Further, $L_{\phi}$ is an $H$-twisted Dirac structure if and only if the spinor $\phi$ satisfies the following condition \cite[Theorem 2.9]{gualtieri:2011a}: 
there exists $A\in \Gamma^\infty(\field TM)$ such that
\begin{equation*}
    \D_{H} \phi=\rho(A)\phi, \quad \textrm{where}\ \D_H:=\D-H\wedge.
 \end{equation*}

\begin{example}\label{example:spinors}
\begin{enumerate}
    \item A spinor of $TM\subseteq \field T M$ is given by the constant function $1\in \Omega^\bullet(M)$, and the spinor line is $\underline{\mathbb{R}}\subseteq \wedge^{\bullet}T^*M$. 
    \item A local spinor of $T^\ast M \subseteq \field T M$ is given by any local volume form on $M$, and the spinor line is $\wedge^{\mathrm{top}}T^*M\subseteq \wedge^{\bullet}T^*M$.
    \item  Let $ (M,\pi) $ be twisted Poisson. On any orientable open subset of $ M $ with volume form $\lambda $ a spinor for $ \graph(\pi) $ is given by 
\begin{equation*}
		\phi_\pi= \E^{\mathrm{i}_\pi} \lambda=\sum_{k=0}^\infty \frac{1}{k!}\I_{\pi}^k \lambda.
\end{equation*}
Here we adopt the convention $\I_{X\wedge Y}=\I_Y \I_X$ for the insertion of multivector fields. 
Further, any  spinor for $ \graph(\pi)$ over such an open subset has this form, for some volume form $\lambda$. 
In particular, the top-degree component of the spinor is nowhere vanishing.
\item Let $\omega\in \Omega^2(M)$ be a $2$-form. A spinor for $\graph(\omega)$ is given by
\begin{equation*}
    \phi_\omega=\E^{\omega}=\sum_{k=0}^\infty \frac{\omega^{\wedge k}}{k!}.
\end{equation*}

\end{enumerate}
\end{example}

We give a useful criterion for extending Dirac structures (for similar results, see \cite{blohmann:2017a}).

\begin{lemma}\label{lem:extending_Dirac_using_spinors}
Let $M$ be a smooth manifold, $H$ a closed 3-form on $M$, and $U\subseteq M$ an open and dense subset. Let $L$ be an $H\at{U}$-twisted Dirac structure on $U$, with spinor line $\Sigma\subseteq \wedge^{\bullet} T^*U$. The following are equivalent.
\begin{itemize}
    \item $L$ extends to an $H$-twisted Dirac structure $\tilde{L}$ on $M$;
    \item $\Sigma$ extends to a line bundle $\tilde{\Sigma}$ on $M$.
\end{itemize}
\end{lemma}

\begin{proof}
Clearly, if $L$ extends, then so does $\Sigma$. Conversely, assume that $\Sigma$ extends to $\tilde{\Sigma}$. We claim that any (local) non-vanishing section $\phi$ of $\tilde{\Sigma}$ is a pure spinor. Indeed, $ L_\phi $ is the kernel of a vector bundle map, so its dimension is locally non-increasing, and, at the same time, $ L_\phi $ is almost everywhere maximally isotropic. This shows that the annihilator $\tilde{L}$ of $\tilde{\Sigma}$ is a Lagrangian extension of the twisted Dirac structure $L$ to $M$. Now, $\tilde{L}$ is also twisted Dirac, because this condition is closed.
\end{proof}

\subsection{Splitting theorems for Dirac structures}
We briefly discuss invariant and transverse submanifolds of Dirac structures, and their normal forms, which play a significant role throughout the paper.

\begin{definition}\label{def:invtrans}
    Let $L\subseteq \field T M$ be a Dirac structure and $N\subseteq M$ a closed and embedded submanifold. 
   \begin{itemize}
               \item $N$ is called \emph{invariant} if $\pr_{TM}(L)$ is tangent to $N$, i.e.\,for all $q\in N$
        \begin{equation*}
            \pr_{TM}(L_q)\subseteq T_qN.
        \end{equation*}
        \item $N$ is called \emph{transversal} if $\pr_{TM}(L)$ is tranverse to $N$, i.e.\ for all $q\in N$
        \begin{equation*}
            T_qN+\pr_{TM}(L_q)=T_qM.
        \end{equation*}
   \end{itemize}
\end{definition}

 Geometrically, invariant submanifolds are those that are unions of presymplectic leaves and transverse ones are those that intersect leaves transversally.

To find local spinors we use Blohmann's normal form theorem developed in \cite[Corollary 3.9]{blohmann:2017a}.  Recall that isomorphisms of twisted Dirac structure include gauge transformations by 2-forms, which are not necessarily closed (Definition \ref{def:isomorphic_dirac_structures}). Therefore, any Dirac structure can be locally untwisted, and so,
Blohmann's splitting theorem can be restated as follows (see also \cite[Theorem 5.1]{bursztyn.lima.meinrenken:2019a}).

\begin{theorem}[\cite{blohmann:2017a}]\label{theorem:blohmann_normal}
Let $L$ be an $H$-twisted Dirac structure on $M$ and let $q\in M$. Then $(M,L)$ is locally isomorphic around $q$ to the (untwisted) product $( U,\graph(\pi))\times (Z,TZ)$ around $(q_{U},q_{Z})$, where $\pi$ a Poisson bivector on $U$ that vanishes at $q_U$.
\end{theorem}

The following observation will be used later on. 

\begin{corollary}\label{corollary:normal_form_invariant}
In the setting of Theorem \ref{theorem:blohmann_normal}, let $N\subseteq M$ be an invariant submanifold through $q$. Then we can take $U=X\times Y$, with $q_U=(q_X,q_Y)$, so that $N$ corresponds to $\{q_X\}\times Y\times Z$, where $\{q_X\}\times Y$ is an invariant submanifold of $U$.
\end{corollary}
\begin{proof}
If $Z$ is connected, invariant submanifolds of $\graph(\pi)\times TZ$ are of the form $Y\times Z$, where $Y$ is an invariant submanifold of $U$. After shrinking $U$, we can choose a complement $X$ to $Y$. 
\end{proof}

\subsection{Invariance and locality of liftability}\label{subs:untwisting}
In the proof of the Main Theorem, we will frequently rely on the local description of twisted Dirac structures, as recalled in the previous subsection. This description is valid up to isomorphisms of Dirac structures (Definition \ref{def:isomorphic_dirac_structures}). In this subsection, we aim to justify the use of this approach. On the one hand, liftability is a local property that is preserved under isomorphisms; on the other hand, the assumptions of the Main Theorem are local in nature and remain unaffected by isomorphisms. We formulate these observations as two lemmas.

\begin{lemma}[Invariance]
For $i=1,2$, let $L_i$ be an $H_i$-twisted Dirac structure on $M_i$. Assume that the two Dirac structures are isomorphic. Fix an isomorphism $(f,B_2)$ consisting of a diffeomorphism $f\colon M_1\diffto M_2$ and a 2-form $B_2\in\Omega^2(M_2)$ such that 
    \begin{equation*}
        L_1=f^\ast(\E^{B_2}L_2) \quad\text{ and }\quad H_1=f^\ast(H_2+\D B_2).
    \end{equation*}
Let $N_1\subseteq M_1$ be a closed and embedded submanifold, and denote $N_2:=f(N_1)$. The following hold.
\begin{enumerate}
        \item $L_1$ lifts to $\blup(M_1,N_1)$ if and only if $L_2$ lifts to $\blup(M_2,N_2)$;
        \item $N_1$ is a transversal for $L_1$ if and only if $N_2$ is a transversal for $L_2$;
        \item $N_1$ is an invariant submanifold for $L_1$ if and only if $N_2$ is an invariant submanifold for $L_2$. Moreover, in this case $f$ induces an isomorphism of bundle of Lie algebras: 
\begin{equation}\label{EQ:Lie_alg_map}(Tf^{-1})^*\colon(TN_1)^{\circ}\diffto (TN_2)^{\circ}.
        \end{equation}
    \end{enumerate}
\end{lemma}
\begin{proof}
For property 1., note that $f$ lifts to a diffeomorphism $\hat{f}\colon\blup(M_1,N_1)\diffto \blup(M_2,N_2)$. Let $p_i\colon\blup(M_i,N_i)\to M_i$ denote the blowdown map. Assume that $L_1$ lifts to a $p_1^*H_1$-twisted Dirac structure $\tilde{L}_1$ on $\blup(M_1,N_1)$. This can be pushed forward via $\hat{f}$ to the Dirac structure $\hat{f}_*(\tilde{L}_1)$ on $\blup(M_2,N_2)$, which is twisted by
\[\hat{f}_*(p_1^*H_1)=\hat{f}_*p_1^*f^*(H_2+\D B_2)=p_2^*(H_2+\D B_2).\]
Then $\tilde{L}_2:=e^{-p_2^*B_2}\circ \hat{f}_*(\tilde{L}_1)$ is a $p_2^*H_2$-twisted Dirac structure which lifts $L_2$. The symmetry of the statement yields the other implication.

For the other properties, note that $(f,B_2)$ induces an isomorphism of Lie algebroids
\[e^{-B_2}\circ \mathbb{T}f\colon L_1\diffto L_2,\qquad
\begin{pmatrix}v \\\xi\end{pmatrix}
\mapsto \begin{pmatrix} Tf(v)\\
(Tf^{-1})^*(\xi)-\I_{ Tf(v)}B_2\end{pmatrix},\]
and that Lie algebroid isomorphisms send transverse (respectively invariant) submanifolds to transverse  (respectively invariant) submanifolds.  Finally, if $N_1$ is invariant, then $(TN_1)^{\circ}\subseteq L_1$ is a subalgebroid with zero anchor (which induces the structure of bundle of Lie algebras), and the restriction of the Lie algebroid isomorphism to this subbundle is precisely \eqref{EQ:Lie_alg_map}.
\end{proof}

The properties above are also local in nature. We state this, but omit the 
obvious proof. 

\begin{lemma}[Locality]
    Let $L$ be an $H$-twisted Dirac structure on $M$ and $N\subseteq M$ a closed and embedded submanifold. The following hold.

\begin{enumerate}
    \item $L$ lifts to $\blup(M,N)$ if and only if every point in $N$ has an open neighbourhood $U$ in $M$ such that $L\at{U}$ lifts to $\blup(U,N\cap U)$;
\item $N$ is a transversal for $L$ if and only if every point in $N$ has an open neighbourhood $U$ in $M$ such that $N\cap U$ is a transversal for $L\at{U}$;
   \item $N$ is an invariant submanifold for $L$ if and only if every point in $N$ has an open neighbourhood $U$ in $M$ such that $N\cap U$ is an invariant submanifold for $L\at{U}$. 
\end{enumerate}
\end{lemma}

\section{Blowup of transversals}\label{sec:blowup_trans}
In this section, we show that any Dirac structure lifts to the blowup along a transversal.

We will use the following terminology. A smooth map $f\colon B\to M$ is said to be \emph{transverse} to a Dirac structure $L$ on $M$, if, for all $b\in B$, the following holds:
\[\pr_{TM}(L_{f(b)})+\image T_b f =T_{f(b)} M.\]
Geometrically, this means that $f$ is transverse to all leaves of $L$. 

The following characterisation and property of transverse maps are well-known (see e.g.\ \cite[Lemma 1.9]{alekseev.bursztyn.meinrenken:2009a} and \cite[Proposition 5.6]{bursztyn:2013a}). For the convenience of the reader, a short proof is included.

\begin{lemma}\label{lemma:transvers_spinor}
    Let $L$ be a twisted Dirac structure on $M$, and $f\colon B\to M$ a smooth map. The following are equivalent:
    \begin{itemize}
\item $f$ is transverse to $L$;
\item for any (local) spinor $\phi$ defining $L$, $f^*\phi$ is nowhere zero.
    \end{itemize}
In this case, the pullback $f^!L$ is a smooth Dirac structure defined by the (local) spinor $f^*\phi$. 
\end{lemma}

\begin{proof}For $b\in B$ denote $V_b:=\pr_{TM}(L_{f(b)})$. Any spinor for $L_{f(b)}$ can be written as $\phi_{f(b)}=\E^{B}\theta_1\wedge\ldots \wedge \theta_n$, where $B\in \wedge^2T^*_{f(b)}M$ and $\theta_1,\ldots,\theta_n$ are a basis of $V_b^{\circ}$, the annihilator of $V_b$ \cite[Proposition 1.3]{gualtieri:2011a}. Then 
\[(T_bf)^*\phi_{f(b)}=\E^{(T_bf)^*B}(T_bf)^*(\theta_1)\wedge\ldots \wedge (T_bf)^*(\theta_n)\] is non-zero precisely when $(T_bf)^*\colon V^{\circ}\to T^*_{b}B$ is injective, which is equivalent to the map $T_bB\to T_{f(b)}M/V_b$ induced by $T_bf$ being surjective, which is equivalent to $V_b+\image T_bf=T_{f(b)}M$. This proves the equivalence.

If $f$ is transverse to $L$, then $L$ can be pulled back along $f$ to a Dirac structure $f^!L$ on $B$ (see e.g.\ \cite[Proposition 5.6]{bursztyn:2013a}). Let $\phi$ be a spinor for $L$. By the first part, $f^*\phi$ is nowhere vanishing. Elements of $f^!L$ have the form $(v,(Tf)^*\eta)$, where $(Tf(v),\eta)\in L$. This implies that 
\[(v,(Tf)^*\eta)\cdot f^*\phi=(Tf)^*((Tf(v),\eta)\cdot\phi)=0.\]
So at every point $f^!L$ is included in the isotropic subspace associated to $f^*\phi$. Since the first vector space is maximally isotropic, the two must be equal. So $f^*\phi$ must be a spinor for $f^!L$ and $f$ is a backward Dirac map.
\end{proof}

Lemma \ref{lemma:transvers_spinor} immediately implies that Dirac structures always lift to the blowup of transversals. 
\begin{theorem}\label{thm:lift_to_transversals}
Let $L$ be a twisted Dirac structure on $M$. For any closed, embedded trans\-versal $N\subseteq M$, $L$ lifts to $\blup(M,N)$ and the blowdown map is a backward Dirac map.  
\end{theorem}
    
\begin{proof}
Since the blowdown map $p$ satisfies
$T_{p(\xi)}N\subseteq \image T_{\xi}p$,  for all $\xi \in\sdiv$, if follows that $p$ is transverse to $L$. The claim follows from Lemma \ref{lemma:transvers_spinor}.
\end{proof}

\section{The transverse and invariant dichotomy}\label{section:dichotomy}

In this section, we prove a surprising result.

\begin{theorem}\label{theorem:Dirac_everything}
    Let $L$ be a twisted Dirac structure on $M$ and $N\subseteq M$ a closed, embedded, and connected submanifold of codimension $>1$. If $L$ lifts to a twisted Dirac structure on $\blup(M,N)$, then one of the following two conditions holds.
    \begin{itemize}
        \item $N\subseteq M$ is a transversal. In this case, the blowdown map is a backward Dirac map.
        \item $N\subseteq M$ is an invariant submanifold.
        In this case, the blowdown map is a forward Dirac map.
    \end{itemize}
\end{theorem}

Motivated by Lemma \ref{lem:extending_Dirac_using_spinors}, we treat the extension problem of line bundles in the next subsection first, and then return to the proof of Theorem \ref{theorem:Dirac_everything}.

\subsection{Extending line bundles along codimension one submanifolds}\label{section:extending_line_bundles}

Aiming to understand the question of extending spinor lines, in this subsection we provide necessary and sufficient conditions for when a rank one subbundle of a vector bundle, defined away from a codimension one submanifold, extends 
smoothly along the submanifold. 

We first make precise the notion of vanishing order. Let $E\to M$ be a vector bundle, and denote by $J^kE\to M$ the \emph{$k$-th order jet bundle} of $E$. This is again a vector bundle over $M$ with fibre at $q\in M$ given by
\[  J^k_qE=\Gamma^\infty(E) / \vanishing_{q}^{k+1} \Gamma^\infty(E),\]
where $\vanishing_q$ the ideal of functions vanishing at $q$. Then any smooth section $s\in \Gamma^{\infty}(E)$ induces a smooth section $j^k(s)\in \Gamma^{\infty}(J^kE)$, called the \emph{$k$-th order jet of $s$}, via  
\[j^k(s)(q):=s+\vanishing_{q}^{k+1} \Gamma^\infty(E)\in J^k_qE.\]

\begin{definition}
The \emph{vanishing order} of a section $s\in \Gamma(E)$ at $q\in M$ is the smallest $k\in \mathbb{N}_0$ such that $j^k(s)(q)\neq 0$. We say that $s$ has \emph{infinite vanishing} order at $q$ if $j^k(s)(q)=0$ for all $k$.
\end{definition}

The following result will be used in Sections \ref{section:proof_of_dirac_everything} and  \ref{section:Poisson_structures_by_spinors}.

\begin{lemma}\label{lemma:extending_line_bundles}
	Let $ E\to M $ be a vector bundle, $ X\subseteq M $ a closed, embedded, and connected submanifold of codimension one. Let $s\in \Gamma^\infty(E)$ be given, satisfying 
    $s^{-1}(0)\subseteq X$.
	\begin{enumerate}
		\item If $s$ has constant finite vanishing order along $X$, then $\Span(s)\at{M\backslash X}$ extends along $X$ to a smooth subbundle of $E$ of rank one.
        \item If $\Span(s)\at{M\backslash X}$ extends along $X$ to a smooth subbundle of $E$ of rank one, and there exists a nowhere vanishing section $\beta \in \Gamma^\infty(E^\ast)$ such that $\beta(s)$ has constant finite vanishing order along $X$, then $s$ has (possibly lower) constant vanishing order along $X$.
	\end{enumerate}
\end{lemma}

For the proof of Lemma \ref{lemma:extending_line_bundles} we need an auxiliary statement about quotients of smooth functions.

\begin{lemma}\label{lemma:quotient_of_smooth_functions} Let $ M $ be a manifold and $X\subseteq M$ a closed, embedded, and connected submanifold of codimension one. Suppose that $f\in \Cinfty(M)$ has constant vanishing order $k\in \field N_0$ along $X$. If $ g\in \Cinfty(M) $ divides $f$, then $g$ has constant vanishing order $\ell\leq k$ along $X$.
\end{lemma}

\begin{proof}
  We can argue locally, and assume that $M=\field R^n$ with coordinates $x_1,\ldots,x_n$ and that $X=\{x_1=0 \}$. Since $f$ vanishes to order $k$ along $X$, by Hadamard's Lemma we can write 
    \begin{equation*}
        f=x_1^{k}\tilde{f}
    \end{equation*}
    for some $\tilde{f}\in \Cinfty(\field R^n)$. Since $j^kf$ vanishes nowhere along $X$, it follows that $\tilde{f}\at{X}$ vanishes nowhere on $X$. Let $U$ be a neighbourhood of $X$ on which $\tilde{f}$ is invertible.  
    
Let $g\in \Cinfty(\field R^n)$ divide $f$. Since $\tilde{f}\at{U}$ is invertible, there is $h\in \Cinfty(U)$ such that 
\[x_1^k=g\cdot h\quad \textrm{on}\quad U.\]

For any $y\in X$ and a function $a\in \Cinfty(U)$, let $o_y(a)$ denote the vanishing order of the function $x_1\mapsto a(x_1,y)$ at $x_1=0$; in other words $o_y(a)$ is the degree of the first Taylor term with non-zero coefficient. Then $o_y$ is multiplicative, $o_y(ab)=o_y(a)+o_y(b)$. Therefore, we obtain that:
\[\forall \ y\in X\ : \ o_y(g)+o_y(h)=k.\]
On the other hand, note that for any $a\in \Cinfty(U)$ the function $y\mapsto o_y(a)$ is locally non-increasing.

Since $X$ is connected, this implies that $o_y(g)$ and $o_y(h)$ are constant. So there is $0\leq \ell\leq k$ such that $o_y(g)=\ell$, for all $y\in X$. 
Thus, we can write $g=x_1^{\ell}\tilde{g}$, where $\tilde{g}\at{X}$ is a nowhere vanishing function. Hence, $g$ has constant vanishing order $\ell$ along $X$. 
\end{proof}

\begin{proof}[of Lemma \ref{lemma:extending_line_bundles}]
	First note that if an extension of $ \Span (s) $ to $X$ exists, it is necessarily unique since $ M\setminus X $ is dense. For the first part, assume that $s$ has constant vanishing order $k$ along $X$. 
    Then for any local defining function $ x $ of $ X $, we have that
	\begin{equation*}
	\tilde{s}=\tfrac{s}{x^k}.
	\end{equation*}
	extends along $X$ to a smooth section of $ E $, and does not vanish anywhere in a neighbourhood of $X$.
    Hence, $\tilde{s}$ generates a line bundle that gives the desired extension. 
    
    For the second part, suppose $ \tilde{s}\in \Gamma^\infty(E) $ is a (local) frame of the extension of $\Span(s)\at{M\backslash X}$. Then there exists a function $ g\in \Cinfty(M) $ with 
	\begin{equation*}
	s=g\tilde{s}.
	\end{equation*}
	Hence $\beta(s)=g\beta(\tilde{s})$, and so $g$ divides $\beta(s)$. If $\beta(s)$ has constant vanishing order $k$ along $X$, by Lemma~\ref{lemma:quotient_of_smooth_functions} $g$ has constant vanishing order $\ell$ along $X$ for some $\ell\leq k$. Then, from $ s=g\tilde{s} $ and since $ \tilde{s} $ does not vanish on $ X$ by assumption, the statement follows.
\end{proof}

The following consequence of our techniques will be used to show that Dirac structures can be extended to an open and dense subset inside the blowup.
 
\begin{lemma}\label{corollary:extending_line_bundle_to_dense_set}
Let $E\to M $ be a vector bundle and $X\subseteq M$ a closed embedded submanifold of codimension one. Let $ s\in \Gamma^\infty(E) $ be a section such that  $s^{-1}(0)\subseteq X$ and $ j^\infty(s)$ is nowhere vanishing. There exists an open dense set $V\subseteq X$ such that $\Span(s)\at{M\backslash X}$ extends to a line bundle over $(M\setminus X)\cup V$.
\end{lemma}

\begin{proof}
	From the proof of Lemma \ref{lemma:extending_line_bundles} we see that, if $s$ has constant vanishing order $\ell\geq 0$ along an open subset $V_\ell\subseteq X$, then $\Span(s)\at{M\backslash X}$ extends to $(M\backslash X)\cup V_\ell $. In particular, we can choose
	\begin{equation*}
		V_\ell:= C_{\ell}^\interior \setminus C_{\ell+1},\quad \textrm{where} \quad C_0:=X,\quad   C_{\ell+1}:=j^{\ell} s\at{X}^{-1}(0), 
	\end{equation*}
and $\argument^{\circ}$ the interior  as a subset of $X$.   
	Hence, if we set $ V=\bigcup_{\ell=0}^\infty V_\ell $, then $\Span(s)\at{M\backslash X}$ extends to a smooth line bundle on $(M\setminus X)\cup V $.
    
    It is left to show that $V$ is dense in $X$. First, using that $ C_{\ell+1}\subseteq C_{\ell}$, that $C_0=X$, and that, by assumption, $ \bigcap_{\ell=0}^{\infty}C_{\ell} =\emptyset$, we obtain the following decomposition of
    $X$ into disjoint subsets:
    \[X=\bigcup_{\ell\geq 0} C_\ell\backslash C_{\ell+1}.\]
    Using that $V_{\ell}\subseteq C_{\ell}$, we have that
    \[X\backslash V=\bigcup_{\ell\geq 0} \big (C_\ell\backslash C_{\ell+1}\big)\backslash V_{\ell}=
    \bigcup_{\ell\geq 0} C_\ell\backslash (C_{\ell}^{\circ}\cup C_{\ell+1})=
    \bigcup_{\ell\geq 0} \partial C_\ell\backslash C_{\ell+1}.\]
	Since all $C_{\ell}$ are closed, their boundaries are closed sets with no inner points. By Baire's Theorem, their union still has no inner points. The same holds also for $\cup_{\ell\geq 0} \partial C_\ell\backslash C_{\ell+1}$, and so, $ V $ is dense.
\end{proof}

Lemma \ref{corollary:extending_line_bundle_to_dense_set} implies the following result. 
   \begin{corollary}\label{coro:almost_everywhere}
Let $L$ be a twisted Dirac structure on $M$ and $N\subseteq M$ a closed and embedded submanifold.  There exists an open and dense subset
   $V\subseteq \field P(\nu_N(M))$ such that the Dirac structure $p^\ast(L\at{M\setminus N})$
    extends smoothly to a Dirac structure on $V\cup \big(\blup(M,N)\setminus \field P(\nu_N(M))\big)$.
\end{corollary}
\begin{proof}Let $\Sigma$ be the spinor line on $\big(\blup(M,N)\setminus \field P(\nu_N(M))\big)$ of $p^*(L\at{M\backslash N})$. Using the standard charts on the blowup, in which the blowdown map is algebraic \eqref{eq:blowdowncoord}, it is easy to check that for any (local) spinor $\phi$ for $L$, we have that $j^{\infty}p^*(\phi)$ does not vanish anywhere. The previous lemma provides an open and dense subset $V\subseteq \field P(\nu_N(M))$ over which $\Sigma$ extends smoothly. Lemma \ref{lem:extending_Dirac_using_spinors} implies that $p^*(L\at{M\backslash N})$ extends to a $p^*H$-twisted Dirac structure to $V\cup \big(\blup(M,N)\setminus \field P(\nu_N(M))\big) $. 
\end{proof}

\subsection{Proof of Theorem \ref{theorem:Dirac_everything}}\label{section:proof_of_dirac_everything}

We first prove a pointwise version of Theorem \ref{theorem:Dirac_everything}. 

\begin{lemma}\label{lemma:coreg_point}
    Let $L$ be a Dirac structure on $M$ and $N\subseteq M$ a closed, embedded submanifold of codimension $>1$. If $L$ lifts to $\blup(M,N)$, then at any $q\in N$ either of the following conditions holds.
    \begin{itemize}
\item\label{lemma:coreg_point:item:inv} $\pr_{TM}(L_q)\subseteq T_q N$;
        \item\label{lemma:coreg_point:item:trans} $\pr_{TM}(L_q)+T_qN=T_q M$.
     \end{itemize}
\end{lemma}
\begin{proof} As explained in Subsection \ref{subs:untwisting}, it suffices to work locally, with an untwisted Dirac structure.

Assume that $L$ lifts, but neither condition holds at $q\in N$. The failure of the first condition implies the existence of some $v+\alpha\in \Gamma^\infty(L)$ such that $v_q\notin T_q N$.

We follow the first steps in \cite{blohmann:2017a} of the proof of Theorem \ref{theorem:blohmann_normal}. By \cite[Lemma 3.4]{blohmann:2017a}, there exists a closed $2$-form $\omega$ defined around $q\in M$ such that $\I_v\omega=\alpha$. Consequently, around $q$ the Dirac structure $L$ is isomorphic to $\E^{-\omega}L$, which contains $v+0$. Hence, we can assume $v\in \Gamma^\infty(L)$ from the start. 
    By using a chart around $q$ adapted to $N$, we find a small neighbourhood $U$ of $q$ in $M$ and a codimension-one submanifold $\iota\colon \tilde{M}\hookrightarrow M$ through $q$, such that $N\cap U\subseteq \tilde{M}$ and such that $v_q\notin T_q\tilde{M}$. By \cite[Lemma 3.5]{blohmann:2017a}, after shrinking $U$ 
    and $\tilde{M}$, we find a diffeomorphism $U\cong (-\varepsilon, \varepsilon)\times \tilde{M}$, which sends $\tilde{M}$ to $\{0\}\times \tilde{M}$, $v\at{U}$ to $\frac{\partial}{\partial t}$, and the Dirac structure $L\at{U}$ to the product Dirac structure 
    \begin{equation*}
       T(-\varepsilon,\varepsilon)\times \tilde{L},
       \end{equation*}
    where $\tilde{L}=\iota^! L$. We replace $M$ by this neighbourhood. Then, because of the product structure, the spinor $\phi_L$ of $L$ may be taken to be the spinor $\phi_{\tilde{L}}$ of $\tilde{L}$. In particular, $\phi_L$ is independent of $t$ and $\D t$. 

Since the second condition
does not hold, and $(\pr_{TM}(L_q))^\circ=L\cap T_q^\ast M$, there exists $\xi\neq 0$, such that $\xi\in L_q\cap (T_qN)^\circ$. Since $\xi,\frac{\partial}{\partial t}\at{q}\in L$, we have that $\xi(\frac{\partial}{\partial t})=0$, and therefore $\xi\at{T_{q}\tilde{M}}\neq 0$. 
So there exists a function $x_2\in \Cinfty(\tilde{M})$ with $x_2\at{N}=0$ and $\xi=\D x_2\at{q}$. By completing $\{x_2\}$ to a chart on $\tilde{M}$ which is adapted to $N$ and centred at $q$, we obtain coordinates $(t,x_2,\dots,x_k,y_{k+1},\ldots,y_{m})$ in which $N=\{t=x_2=\dots=x_k=0\}$ and $t(q)=x_i(q)=y_j(q)=0$.
    
In the constructed chart, the spinor $\phi_L$ of $L$ decomposes uniquely as
    \begin{equation*}
        \phi_L=\phi_0+\phi_1\wedge \D x_2,
    \end{equation*}
where $\phi_0$ has no $\D x_2$-contribution. We still have that $\phi_0,\phi_1$ are independent of $t$ and $\D t$, and,  
since $\D x_2\at{q}\in L$, by the definition of the action $\rho$ in Section \ref{sec:spinor}, we have $\phi_0(0)=0$.  Next, we use multi-indices to write
    \begin{equation*}
        \phi_0=\sum_{I,J}\phi_0^{I,J} \D x_I\wedge \D y_J\quad \text{ and } \quad \phi_1=\sum_{I,J}\phi_1^{I,J} \D x_I\wedge \D y_J.
    \end{equation*}
We pull the spinor $\phi_L$ back to the chart $U_t$ on $\blup(M,N)$, i.e.\ along the map $p(\tilde{t},\tilde{x},y)=(\tilde{t},\tilde{t}\tilde{x},y)$, 
where $\tilde{x}=(\tilde{x}_2,\ldots,\tilde{x}_k)$ and $y=(y_{k+1},\ldots,y_m)$ (see \eqref{eq:blowdowncoord}).
Then we obtain 
    \begin{equation*}
        \begin{aligned}           p^\ast\phi_L\at{U_t}=&\sum_{I,J}\Big( \tilde{t}^{|I|}p^\ast\phi_0^{I,J} \D \tilde{x}_I\wedge \D y_J + \tilde{t}^{|I|+1}p^\ast\phi_1^{I,J}\D \tilde{x}_I\wedge \D y_J\wedge \D \tilde{x}_2 \\
            &+\D \tilde{t}\wedge \Big( \sum_{u=3}^k \tilde{t}^{|I|-1}p^\ast\phi^{I,J}_0 \tilde{x}_u\I_{\partial_{\tilde{x}_u}} \D \tilde{x}_I \wedge \D y_J +\sum_{u=3}^k  \tilde{t}^{|I|}p^\ast\phi^{I,J}_1 \tilde{x}_u\I_{\partial_{\tilde{x}_u}} \D \tilde{x}_I \wedge \D y_J \wedge \D \tilde{x}_2\\
            &+(-1)^{|I|+|J|}\tilde{x}_2 \tilde{t}^{|I|}p^\ast\phi_1^{I,J}\D \tilde{x}_I\wedge \D y_J\Big)\Big).
        \end{aligned}
    \end{equation*}
    Recall that $U_t\cap \field P(\nu_N(M)))=\{\tilde{t}=0\}$.
    Since by assumption the Dirac structure lifts, there exists a function $f\in \Cinfty(U_t\setminus \field \{\tilde{t}=0\})$ such that $f p^\ast\phi_L\at{U_t}$ extends smoothly along $\{\tilde{t}=0\}$, and vanishes nowhere on $U_t$. In particular, for $(\tilde{x},y)=0$ and $\tilde{t}\neq 0$, using that $\phi_0(0)=0$, we see that all terms in of $ p^\ast\phi_L\at{U_t}$ in the sum above vanish identically on the line $(\tilde{x},y)=0$, except for \begin{equation*}
\sum_{I,J}\tilde{t}^{|I|+1}p^\ast\phi_1^{I,J}(0)\D \tilde{x}_I\wedge \D y_J\wedge \D \tilde{x}_2.
    \end{equation*}
   Hence, there must exist index sets $I_\ast, J_\ast$ such that the smooth function 
\[g:=f\tilde{t}^{|I_\ast|+1}p^\ast\phi_1^{I_\ast,J_\ast}\in \Cinfty(U_t)\]
satisfies $g(0)\neq 0$. By shrinking $U_t$, we may assume that $g$ doesn't vanish anywhere on $U_t$, and by replacing $f$ by $f/g$, we may assume that $f\tilde{t}^{|I_\ast|+1}p^\ast\phi_1^{I_\ast,J_\ast}=1$ on $U_t$. Hence,  
    \begin{equation*}
        f=\frac{1}{\tilde{t}^{|I_\ast|+1}p^\ast\phi_1^{I_\ast,J_\ast}}.
    \end{equation*}
    However, this form of $f$ implies that $f p^\ast\phi_L\at{U_t}$ has a singularity at $\tilde{t}=0$. Indeed, for $\tilde{t}\neq 0$, along the submanifold $\{x_3=\ldots=x_k=0\}$, we have 
    \begin{equation*}
        \begin{aligned}
            f p^\ast\phi_L=& \frac{1}{\tilde{t}^{|I_\ast|+1}p^\ast\phi_1^{I_\ast,J_\ast}}\sum_{I,J}\Big(\tilde{t}^{|I|}p^\ast\phi_0^{I,J} \D \tilde{x}_I\wedge \D y_J + \tilde{t}^{|I|+1}p^\ast\phi_1^{I,J}\D \tilde{x}_I\wedge \D y_J\wedge \D \tilde{x}_2 \\
            &+(-1)^{|I|+|J|}\tilde{x}_2 \tilde{t}^{|I|}p^\ast\phi_1^{I,J}\D \tilde{t}\wedge \D \tilde{x}_I\wedge \D y_J\Big).
        \end{aligned}
    \end{equation*}
    The coefficient of the term $\D \tilde{t}\wedge \D \tilde{x}_{I_\ast}\wedge \D y_{J_\ast}$ is given by 
    \begin{equation*}
        \frac{\tilde{x}_2}{\tilde{t}}\D \tilde{t}\wedge \D \tilde{x}_{I_\ast}\wedge \D y_{J_\ast},
    \end{equation*}
    which is not well-defined in $\tilde{t}=0$. Hence, we obtain a contradiction to the assumption that the Dirac structure lifts. 
\end{proof}

Next, we use Lemma \ref{lemma:extending_line_bundles} to show that, if $L$ lifts over a point where the Dirac structure is tangent to the submanifold, then locally the pulled back spinor has constant vanishing order.

\begin{lemma}\label{lemma:invariant_points_are_open}
    Let $L$ be a Dirac structure on $M$, $N\subseteq M$ a closed and embedded submanifold, and $q\in N$ such that 
    \begin{equation}\label{eq:pointwise_tangent}
        \pr_{TM}(L_q)\subseteq T_qN.
    \end{equation}
    Let $\phi$ denote a spinor corresponding to $L$ defined around $q$. The following are equivalent. 
    \begin{itemize}
        \item The Dirac structure $p^\ast(L\at{M\setminus N})$ extends to a neighbourhood of $p^{-1}(q)\subseteq \sdiv$.
    \item The order of vanishing of $p^\ast\phi$ along $\sdiv$ is constant around $p^{-1}(q)$.
    \end{itemize} 
\end{lemma}

\begin{proof} If $p^*\phi$ has constant vanishing order along $\field P(\nu_N(M))$, then the first part of Lemma \ref{lemma:extending_line_bundles} together with Lemma \ref{lem:extending_Dirac_using_spinors} imply that $L$ lifts in a neighbourhood of $p^{-1}(q)$. 

For the converse, we apply the second part of Lemma \ref{lemma:extending_line_bundles}. In order to do so, we first find coordinates that are adapted to both $L$ and $N$.

    By Theorem \ref{theorem:blohmann_normal} (see also Subsection \ref{subs:untwisting}) we can assume that $M=U\times Z$, $q=(q_U,q_Z)$, and $L=\graph(\pi)\times T Z$ for a Poisson structure $\pi$ on $U$ vanishing at $q_U$. 
    
    The assumption \eqref{eq:pointwise_tangent} implies that $\pr_Z\at{N}\colon N\to Z $ is a submersion at $q$. After shrinking $U$, one can split $U=V\times W$, with $q_U=(q_V,q_W)$, such that $T_{q_W}W=\ker T_q (\pr_Z\at{N})$. Then the projection $\pr \colon V\times W\times Z\to W\times Z$ induces a linear isomorphism 
    \begin{equation*}
        T_q\pr\at{N}\colon T_q N\diffto  T_{q_W}W\oplus T_{q_Z} Z.
    \end{equation*}
    By the inverse function theorem, $\pr\at{N}\colon N\to W\times Z$ is a diffeomorphism in an open neighbourhood of $q$. Therefore, after shrinking all neighbourhoods, we find a smooth map $f\colon  W\times Z\to V$, such that
    \begin{equation*}
        (\pr\at{N})^{-1}(w,z)=(f(w,z), w, z).
    \end{equation*}
    
    Denote by $v=(v_1,\ldots, v_k)$, $w=(w_1,\ldots, w_l)$, and $z=(z_1,\ldots, z_m)$ coordinates on $V$, $W$, and $Z$, around $q_V$, $q_W$, and $q_Z$, respectively. Define $x_i:=v_i-f_i(w,z)$, $i=1,\ldots, k$, and $x=(x_1,\ldots,x_k)$. 
    Then $(x,w,z)$ is a chart around $q$ on $M$ which is a submanifold chart for $N$, i.e.\ $N=\{x=0\}$. 
    
    In the first set of coordinates, a spinor of $L$ is given by
    \begin{equation*}
        \phi=\E^{\I_\pi}(\D v_1\wedge\ldots\wedge \D v_k\wedge \D w_1\wedge\ldots\wedge  \D w_l)= 
        \D v_1\wedge\ldots\wedge \D v_k\wedge \D w_1\wedge\ldots\wedge  \D w_l+ \text{ lower degree forms}.
    \end{equation*}
Writing $\phi$ in the second set of coordinates, we obtain
    \begin{equation*}
        \phi=\D x_1\wedge\ldots\wedge \D x_k\wedge \D w_1\wedge\ldots\wedge  \D w_l + \text{ other terms},
    \end{equation*}
where the ``other terms'' are either of lower form degree or contain at least one $\D z_j$. 

For $i=1,\ldots,k$ consider the chart $U_i$ on $\sdiv$ with coordinates $(\tilde{x},w,z)$, defined in \eqref{eq:blowdowncoord}. Then we have that
    \begin{equation}\label{eq:some_pullback_spinor}
p^\ast\phi\at{U_{x_i}}=\tilde{x}_i^{k-1} \D \tilde{x}_1\wedge\ldots\wedge \D \tilde{x}_k\wedge \D w_1\wedge\ldots\wedge  \D w_l + \text{ other terms},
    \end{equation}
    where again the ``other terms'' are either of lower degree or contain at least one $\D z_j$. Hence, the multi-vector field 
    \[\beta=\frac{\partial }{\partial\tilde{x}_1 }\wedge\ldots\wedge \frac{\partial }{\partial\tilde{x}_k }\wedge \frac{\partial }{\partial w_1 }\wedge\ldots\wedge  \frac{\partial }{\partial w_l},\]
    satisfies that $\beta(p^*\phi)$ has constant vanishing order equal to $k-1$ along $\sdiv\cap U_i=\{\tilde{x}_i=0\}$. So if the line bundle spanned by $p^*\phi$ extends to a neighbourhood of $p^{-1}(q)$, then the second part of Lemma \ref{lemma:extending_line_bundles} implies that $p^*\phi$ has constant vanishing order in a neighbourhood of $p^{-1}(q)$.  
\end{proof}

The following is used to show that the blowdown map is a forward Dirac map. 

\begin{lemma}\label{lemma:general_dirac:forward_dirac_equivalent_to_invariance}
         Let $L$ be a Dirac structure on $M$ and $N\subseteq M$ a closed, embedded submanifold of codimension $>1$. Assume that $L$ lifts to the Dirac structure $\tilde{L}$ on $\blup(M,N)$. Then the blowdown map is a forward Dirac map if and only if $N$ is an invariant submanifold.
    \end{lemma}

    \begin{proof}
        Suppose $\blowdown{p}{M}{N}$ is a forward Dirac map.
        For any $ q\in N $ we have 
        \begin{equation*}
            \pr_{TM}(L_q)\subseteq \bigcap_{ \xi\in p^{-1}(q) }\image T_{\xi}p=T_qN.
        \end{equation*} 
        Hence, $N$ is invariant. 
        
        Conversely, assume that $N$ is invariant.
        We have to show that for all $ \xi\in \field P(\nu_{N}(M)) $, $p_!(\tilde{L}_{\xi})=L_{p(\xi)}$. 
Fix $(z,\eta)\in L_{p(\xi)}$. Let $(X,\alpha)\in \Gamma^\infty(L)$ with $(X,\alpha)\at{p(\xi)}=(z,\eta)$. Since $X$ is tangent to $N$, by Lemma \ref{lemma:lifting_vf_to_the_blowup}, it has a unique lift $\tilde{X}$ to $\blup(M,N)$. 
Then $(\tilde{X},p^\ast\alpha)\in \Gamma^\infty(\field T\blup(M,N))$ maps the dense subset $\blup(M,N)\setminus \sdiv$ into $\tilde{L}$. 
Thus, $(\tilde{X},p^\ast\alpha)\in \Gamma^\infty(\tilde{L})$, implying that
\begin{equation*}
     \begin{pmatrix}(T_{\xi}p)(\tilde{X}(\xi))\\ \alpha(p(\xi))\end{pmatrix} \in p_!(\tilde{L}_{\xi}). 
\end{equation*}
As $\tilde{X}$ $p$-projects to $X$, we see that the above pair is just $(X,\alpha)\at{p(\xi)}=(z,\eta)$. 
This shows $ p_!(\tilde{L}_{\xi}) \supseteq L_{p(\xi)}$.
Equality follows because both are maximally isotropic subspaces of $\field T_{p(\xi)}M$.
\end{proof}

Summarising, we can prove Theorem \ref{theorem:Dirac_everything}.

\begin{proof}[of Theorem \ref{theorem:Dirac_everything}]
    By Lemma \ref{lemma:coreg_point}, at any $q\in N$, $L$ is either tangent or transverse to $N$. These two cases are mutually exclusive as $\codim N>0 $. Moreover, the set where $N$ is transverse to $L$ is open in $N$, because the condition $T_qN+\pr_{TM}(L_q)=T_qM$ is open. Since $N$ is connected, it suffice to show that the tangential points also form an open subset. 
    Indeed, let $q\in N$ be such that $\pr_{TM}(L_q)\subseteq T_qN$ and let $\phi$ be a spinor for $L$ defined around $q$. By Lemma \ref{lemma:transvers_spinor}, $p^*\phi$ must vanish along $p^{-1}(q)$ as $q$ is not a transverse point. 
    But then, by Lemma \ref{lemma:invariant_points_are_open},  $p^*\phi$ must vanish also in a neighbourhood $U$ of $p^{-1}(q)\subseteq \sdiv$. 
    Then $p(U)\subseteq N$ is a open neighbourhood of $q$ in $N$, and $L$ and $N$ cannot be transverse at any point $q'\in p(U)$, again because of Lemma \ref{lemma:transvers_spinor}. Therefore, Lemma \ref{lemma:coreg_point} shows that $p(U)$ consists entirely of points where $L$ is tangent to $N$. 

Finally, if $N$ is a transversal, Lemma \ref{lemma:transvers_spinor} shows that $p$ is a backward Dirac map, and if $N$ is invariant, then Lemma \ref{lemma:general_dirac:forward_dirac_equivalent_to_invariance} shows that $p$ is a forward Dirac map. 
\end{proof}

\section{Blowup of invariant submanifolds}\label{section:Poisson_structures_by_spinors}

After the previous two sections, it remains to characterise the invariant submanifolds $N\subseteq M$ of a twisted Dirac structure $L$ for which $L$ lifts to $\blup(M,N)$. Recall from the Introduction that an invariant submanifold $N\subseteq M$ of $L$ yields a short exact sequence of Lie algebroids:
\[0\longrightarrow (TN)^{\circ}\longrightarrow L\at{N}\longrightarrow \iota_N^!L\longrightarrow 0.
\]
In particular, $(TN)^\circ$ is a bundle of Lie algebras over $N$. 

Recall also from the Introduction that, given a Lie algebra $(\mathfrak{g},[\argument,\argument])$, the \emph{height} of an element $\xi\in \liealg g^\ast \setminus \{0\}$ is defined as the integer $k\in \field N_0$ such that
\begin{equation}\label{eq:height}
        \xi\wedge (\DEC\xi)^{k}\neq 0 \quad \text{ and }\quad \xi\wedge (\DEC\xi)^{k+1}=0.
    \end{equation}
    The Lie algebra $\liealg g$ is called \emph{a Lie algebra of constant height $k$} if \eqref{eq:height} holds for all $\xi\in \liealg g^\ast\setminus \{0\}$.

We can now state  the main result of this section. 

\begin{theorem}\label{theorem:lifting_by_spinors}  Let $L$ be a twisted Dirac structure on $M$ and $N\subseteq M$ a closed, embedded and connected submanifold, which is invariant. The following are equivalent.
\begin{itemize}
\item $L$ lifts to  $\blup(M,N)$.
\item The Lie algebras $(T_qN)^{\circ}$, $q\in N$, have all the same constant height $k$.
\end{itemize}
\end{theorem}

\begin{example} \label{CartanDirac}
On a Lie group $G$ with a bi-invariant (not necessarily positive) metric  $(\argument,\argument)$, there is a canonical twisted Dirac structure $L$, called Cartan-Dirac structure \cite[Example 5.2]{severa.weinstein:2001a}.
The twist is provided by the corresponding Cartan 3-form, i.e.\ the bi-invariant 3-form on $G$ which at the unit $e$ reads $H(u,v,w):=-\frac{1}{2}([u,v],w )$ for $u,v,w\in \g=T_eG$. In a neighbourhood $U$ of the unit $e$, the Cartan-Dirac structure $L$ is the graph of an $H\at{U}$-twisted Poisson structure $\pi$.
At every $g\in U$, the map $\pi^{\sharp}_g\colon T_g^*G\to T_gG$ is obtained by left-translating  $2(\Ad_g-1)(\Ad_g+1)^{-1}\colon \g\to \g$ and using the identification  $\g^*\cong \g$ induced by the metric. It is known that  the induced Lie algebroid structure on $L$ is isomorphic over $\id_G$ to the transformation Lie algebroid associated to the action of $G$ on itself by conjugation. In particular, the unit element $\{e\}$ of $G$  constitutes a leaf of $L$, and the isotropy Lie algebra of $L$ at $e$ is just $\g$.

Now choose $G=SO(3)$, endowed with any bi-invariant Riemannian metric (it exists because the Lie group  is compact). Theorem \ref{theorem:classification_LA_of_constant_height} below shows that $\mathfrak{so}(3)$ is of constant height 1. Hence, by Theorem \ref{theorem:lifting_by_spinors}, $L$ lifts to an $p^*H$-twisted Dirac structure on $ \blup(G,\{e\})$, which is not given by a bivector field around $p^{-1}(e)$.
\end{example}

Theorem \ref{theorem:lifting_by_spinors} is a consequence of the results of the following subsections, in which we give alternative criteria for the liftability property. 

\subsection{The constant vanishing order criterion}

The following criterion for a Dirac structure to lift to the blowup is an immediate consequence of Lemma \ref{lemma:invariant_points_are_open} and its proof.

\begin{corollary}\label{coro:vanishing:order}
 Let $L$ be a twisted Dirac structure on $M$ and $N\subseteq M$ a closed, embedded and connected submanifold, which is invariant and has codimension $>1$. The following are equivalent. 
    \begin{itemize}
        \item $L$ lifts to a $p^\ast H$-twisted Dirac structure on $\blup(M,N)$.
        \item There exists $\ell\in \{1,\ldots, \mathrm{codim}\, N-1\}$ such that, for any (local) spinor $\phi$ for $L$, $p^*\phi$ has constant vanishing order $\ell$ along $\sdiv$.
        \end{itemize}
\end{corollary}

In the Poisson category, we have the following version of this result. 
\begin{corollary}\label{cor:lift_Poisson_Poisson}
    Let $ (M,\pi) $ be Poisson and $N\subseteq M$ a closed and embedded Poisson submanifold. The following are equivalent.
    \begin{itemize} 
        \item $\pi$ lifts to a  Poisson structure $\tilde{\pi}$ on $\blup(M,N)$.
        \item The spinor $p^*(\E^{\pi}\lambda)$ has constant vanishing order equal to $\ell=\mathrm{codim}\, N-1$ along $\sdiv$, for any (local) volume form $\lambda$ on $M$.
    \end{itemize}
\end{corollary}

\begin{proof}
By Corollary  \ref{coro:vanishing:order}, $\graph(\pi)$ lifts to a Dirac structure $\tilde{L}$, exactly when the spinor $\phi^*(\E^{\pi}\lambda)$ has constant vanishing order $\ell$. In this case, $t^{-\ell}\phi^*(\E^{\pi}\lambda)$ is a nowhere vanishing spinor corresponding to $\tilde{L}$, where $t$ is any local coordinate defining $\sdiv$, i.e.\ locally $\sdiv=\{t=0\}$. On the other hand, $t^{-\ell}\phi^*(\E^{\pi}\lambda)$ corresponds to a Poisson structure, exactly when its component of top degree, i.e.\ $t^{-\ell}\phi^*(\lambda)$, is nowhere vanishing. Using the standard charts on the blowup, as in \eqref{eq:some_pullback_spinor}, we see that $\phi^*(\lambda)$ has constant vanishing order $\codim\,  N -1$. This implies the equivalence.
\end{proof}

\begin{example}\label{example:bundle_of_scaled_so3}
    Consider the vector bundle $E=\field R^3\times \field R^2\to \field R^2$. Denote its canonical global frame by $\{e_1,e_2,e_3\}$ and denote by $(y_1,y_2)$ the variables of the base. Let $f\in \Cinfty(\field R^2)$ be a function. Equip the fibre $E_{(y_1,y_2)}$ with the Lie bracket
    \begin{equation*}
[e_i,e_j]_{(y_1,y_2)}=f(y_1,y_2)\sum_{k=1}^3 \varepsilon_{ijk}e_k,
    \end{equation*}
    where $\varepsilon_{ijk}=1$  if $ijk$ is a cyclic permutation of $123$, $\varepsilon_{ijk}=-1$
 for a cyclic permutation of $213$, and $\varepsilon_{ijk}=0$ otherwise. 
    Then $E\to \field R^2$ is a bundle of Lie algebras, with fibre
    \begin{equation*}
        E_{(y_1,y_2)}\simeq \begin{cases}
            \liealg{so}(3) &\text{ if }f(y_1,y_2) \neq 0\\
            \text{abelian }\field R^3 &\text{ if }f(y_1,y_2)=0.
        \end{cases}
    \end{equation*}
    On $E^\ast$, we obtain a Poisson structure $\pi$ with global spinor line generated by 
    \begin{equation*}
\phi_{\pi}=f(y_1,y_2)\sum_{i=1}^3  x_i\D x_i\wedge \D y_1 \wedge \D y_2 + \D x_1\wedge \D x_2 \wedge \D x_3 \wedge \D y_1 \wedge \D y_2,
    \end{equation*}
    where we denote the linear fibre coordinates induced by (the dual frame to) $\{e_i\}_i$ by $(x_1,x_2,x_3)$.
    \begin{enumerate}
\item\label{example:bundle_of_scaled_so3:item:zero_section} Let $N$ be the zero section of the vector bundle $E^*$.
    The Dirac structure $\graph(\pi)$ lifts to the blowup $\blup(E^\ast,N)$ exactly when $f= 0$ or $f$ vanishes nowhere. Indeed, for $i\in \{1,2,3\}$, 
    \begin{equation*}
        \begin{aligned}
(p^\ast\phi_{\pi})\at{U_{i}} =& f(y_1,y_2)\Big(\sum_{j\neq i}  \tilde{x}_i^2\tilde{x}_j \D \tilde{x}_j +\tilde{x}_i\Big(1+\sum_{j\neq i} \tilde{x}_j^2\Big) \D \tilde{x}_i\Big)\wedge \D y_1 \wedge \D y_2 
            \\&+ \tilde{x}_i^2\D \tilde{x}_1\wedge \D \tilde{x}_2 \wedge \D \tilde{x}_3 \wedge \D y_1 \wedge \D y_2.
        \end{aligned}
    \end{equation*}
Let $q= (\tilde{q}_1,\tilde{q}_2,\tilde{q}_3,z_1,z_2)\in \field P (E^*)\cap U_{i}$, so $\tilde{q}_i=0$, with vanishing ideal $\vanishing_{q}$. Reducing the coefficients in $(p^\ast\phi_{\pi})\at{U_{i}}$ modulo $\vanishing_{q}^2$, the only term that might be non-zero is
\[f(z_1,z_2)\, \tilde{x}_i\, (1+\sum_{j\neq i}\tilde{q}_j^2)\in \vanishing_{q}.\]
This belongs to $\vanishing_{q}^2$ precisely when $f(z_1,z_2)=0$. Hence, the vanishing order of $(p^\ast\phi_{\pi})\at{U_{i}}$ is constant along $\field P (E^*) \cap U_i{=\{\tilde{x}_i=0\}}$ if and only if $f=0$ or $f$ vanishes nowhere. 

If $f=0$, the vanishing order is $2=\mathrm{codim}\, N-1$, $\pi=0$, and the lifted Dirac structure corresponds to the zero Poisson structure. 

If $f$ vanishes nowhere, $(p^\ast\phi_{\pi})\at{U_{i}}$ has vanishing order one along $\sdiv \cap U_i$ and so $\graph(\pi)$ lifts to a Dirac structure, which does not come from a Poisson structure.
    \item Suppose that $f$ vanishes at the origin $0\in {\field R^2}$, and take $N$ to be the fibre $E^\ast\at{0}$. Then, by Corollary \ref{cor:lift_Poisson_Poisson}, $\pi$ lifts to a Poisson structure on $\blup(E^\ast, N)$. Indeed, for $j\in \{1,2\}$,
 \[(p^\ast\phi_{\pi})\at{U_{j}}= \tilde{y}_j \Big(\D x_1\wedge \D x_2 \wedge \D x_3 + \sum_{i=1}^3 (p^\ast f)\at{U_{j}} x_i  \D x_i \Big)\wedge \D \tilde{y}_1 \wedge \D \tilde{y}_2
  \]
  has constant vanishing order $\codim\, N-1=1$ along $\field P(\nu_N(E^*))\cap U_{j}$.
    \end{enumerate}
\end{example}

\subsection{Reducing to bundles of linear Poisson structures}

As a first step to prove Theorem \ref{theorem:lifting_by_spinors}, we reduce the problem to the setting of linear Poisson structures. 

Given a bundle of Lie algebras $(\underline{\mathfrak{g}},[\argument,\argument])$ over $N$, the dual vector bundle $r\colon  \underline{\mathfrak{g}}^*\to N$ comes equipped with a with a fibrewise linear Poisson structure $\plin$, with Poisson bracket determined by: 
 \[\{\ev_{\alpha},\ev_{\beta}\}=\ev_{[\alpha,\beta]},\quad \{\ev_{\alpha},r^*f\}=0, \quad \{r^*f,r^*g\}=0,\]
 for all $\alpha,\beta\in \Gamma^{\infty}(\underline{\mathfrak{g}})$ and all $f,g\in \Cinfty(N)$, where $\ev_{\alpha}\in \Cinfty(\underline{\mathfrak{g}}^*)$ is $\alpha$, viewed as a fibrewise linear map. 

In particular, for an invariant submanifold $N\subseteq M$ of a Dirac structure $L$, the normal bundle
$\nu_N(M)=((TN)^\circ)^*$ has a linear (untwisted) Poisson structure, denoted $\plin$, coming from the bundle of Lie algebras on  $(TN)^\circ$. Linearity implies that $\plin$ vanishes along $0_N\subseteq \nu_N(M)$, and so $0_N$ is an invariant submanifold for $\plin$.

\begin{theorem}\label{pullback of invariant}
    Let $L$ be a twisted Dirac structure on $M$ and $N\subseteq M$ a closed and embedded submanifold,  which is invariant. The following are equivalent. 
 \begin{itemize}
     \item 
  $L$ lifts to $\blup(M,N)$.
   \item $\graph(\plin)$ lifts to $\blup(\nu_N(M),0_N)$.
     \end{itemize}
     In that case, the respective pullback spinors have the same vanishing orders along $\sdiv$. 
\end{theorem}

\begin{proof}
As remarked in Subsection \ref{subs:untwisting}, it suffices to check locally whether $L$ (respectively $\graph(\plin)$) lifts to the blowup. In Subsection \ref{subs:untwisting}, we showed that isomorphisms of Dirac structures induce isomorphisms between the corresponding bundle of Lie algebras $(TN)^{\circ}$, and so, isomorphisms between the corresponding linear Poisson structures $\plin$. Therefore, it suffices to work in the local product neighbourhoods from Corollary \ref{corollary:normal_form_invariant}, as these have isomorphic linear Poisson structures.

So we assume that $M=X\times Y\times Z$, $L=\graph(\pi)\times TZ$, where $\pi$ is a Poisson structure on $X\times Y$, and that $N=\{q_X\}\times Y\times Z$, where $\{q_X\}\times Y\subseteq X\times Y$ is a Poisson submanifold. Moreover, we may assume that $X\subseteq \mathbb{R}^m$, $Y\subseteq \mathbb{R}^n$, 
$Z\subseteq \mathbb{R}^p$ are open subsets, with coordinates $x$, $y$ and $z$, respectively, and that $q_X=0$. Then we can write 
\begin{equation}\label{eq:pi_in_coord}
\pi=\frac{1}{2}\sum_{i,j}\pi_{ij}\frac{\partial}{\partial x_i}\wedge \frac{\partial}{\partial x_j}+ \sum_{i,\alpha} \pi_{i\alpha }\frac{\partial}{\partial x_i}\wedge \frac{\partial}{\partial y_\alpha} +\frac{1}{2}\sum_{\alpha,\beta}\pi_{\alpha\beta}\frac{\partial}{\partial y_\alpha}\wedge \frac{\partial}{\partial y_\beta},
    \end{equation}
where the coefficients are smooth functions of $(x,y)\in X\times Y$. Invariance of $\{0\}\times Y$ means that $\pi_{ij}(0,y)=0$ and $\pi_{i\alpha}(0,y)=0$. We can identify $\nu_N(M)=\mathbb{R}^{m}\times Y\times Z$ and then the Poisson structure $\plin$ becomes the linearisation of the first term: 
\begin{equation*}
    \plin=\frac{1}{2}\sum_{i,j,k}x_k\frac{\partial \pi_{ij}}{\partial x_k}(0,y)\frac{\partial}{\partial x_i}\wedge \frac{\partial}{\partial x_j}.
\end{equation*}
We claim that $w:=\pi-\plin$ can be lifted to a smooth bivector field $\tilde{w}$ on $\blup(M,N)$. This follows from Lemma \ref{lemma:lifting_vf_to_the_blowup} because $w$ can be written as a sum $w=\sum_k U_k\wedge V_k$, where $U_k$ and $V_k$ are vector fields tangent to $N=\{0\}\times Y\times Z$. This is clear for the last term in \eqref{eq:pi_in_coord}, for the middle it follows because $\pi_{i\alpha}(0,y)=0$, and for the remainder it follows because its coefficients vanish quadratically along $N$. 

Consider the differential forms
\[\lambda_X:=\D x_1\wedge \ldots \wedge \D x_m,\quad \lambda_Y:= \D y_1\wedge \ldots\wedge \D y_n, \quad\text{and} \quad  \lambda_Z:=  \D z_1\wedge \ldots\wedge \D z_p.\]
Spinors for $L=\graph(\pi)\times TZ$ and $\graph(\plin)$ are given, respectively, by  
\[\phi_L=\E^\pi(\lambda_{X}\wedge\lambda_Y) = \E^w((\E^{\plin}\lambda_{X})\wedge\lambda_Y)\quad \textrm{and} \quad \phi_{\plin}=(\E^{\plin}\lambda_{X})\wedge\lambda_Y\wedge\lambda_Z.
\]
Their pullbacks to $\mathbb{B}:=\blup(X\times Y\times Z, \{0\}\times Y\times Z)\cong \blup(X,\{0\})\times Y\times Z$ are given by
\begin{align*}
p^*(\phi_L)&=p^*\big(\E^w((\E^{\plin}\lambda_{X})\wedge\lambda_Y\big)=
\E^{\tilde{w}}\big(p^*\big(\E^{\plin}\lambda_{X})\wedge\lambda_Y\big)\\
\textrm{and}\quad    p^*(\phi_{\plin})&=p^*(\E^{\plin}\lambda_{X})\wedge\lambda_Y\wedge\lambda_Z.
\end{align*}

By Corollary  \ref{coro:vanishing:order}, the first condition in the theorem is equivalent to
the vanishing order of $p^*(\phi_L)$ being constant along $\field P(\mathbb{R}^m)\times Y\times Z$. Since $\E^{\tilde{w}}$ is an automorphism of the vector bundle 
$\wedge^{\bullet} T^*\mathbb{B}$, this condition is equivalent to $p^*(\E^{\plin}\lambda_X)\wedge \lambda_Y$ having constant vanishing order. This is then equivalent to $p^*(\phi_{\plin})$ having constant vanishing order, which by Corollary \ref{coro:vanishing:order}, is equivalent to the second condition in the theorem.
\end{proof}

\subsection{Blowup of a bundle of linear Poisson structures}\label{sec:transverse_lie_algebra}

In the previous section, we have reduced the problem of blowing up invariant submanifolds from twisted Dirac structures to bundles of linear Poisson structures, which we will discuss here. The next theorem, combined  with Theorem \ref{pullback of invariant}, yields  Theorem \ref{theorem:lifting_by_spinors}.

\begin{theorem}\label{theorem:summary_submanifold}
Let $(\underline{\mathfrak{g}}^*,\plin)$ be a bundle of linear Poisson structures corresponding to the bundle of Lie algebras $(\underline{\mathfrak{g}},[\argument,\argument])$ over $N$. The following are equivalent.
  \begin{itemize}         \item The Dirac structure $\graph(\plin)$ lifts to $\blup(\underline{\mathfrak{g}}^*,0_N)$.
\item The Lie algebras $\underline{\mathfrak{g}}\at{q}$, $q\in N$, have all the same constant height $k$.
\end{itemize}
\end{theorem}

The first step in the proof is to relate the vanishing orders of $\plin$ to those of its restriction to the fibres. 

\begin{lemma}\label{lemma:normal_jet_to_fibre_jet}
   Let $(\underline{\mathfrak{g}}^*,\plin)$ be a bundle of linear Poisson structures corresponding to the bundle of Lie algebras $(\underline{\mathfrak{g}},[\argument,\argument])$ over $N$. Let $\ell \geq 0$. The following are equivalent.
      \begin{itemize}
        \item $p^*(\E^{\plin}\lambda)$ has constant vanishing order $\ell$ along $\field P(\underline{\mathfrak{g}}^*)$, for any (local) volume form $\lambda$ on $\underline{\mathfrak{g}}^*$.
        \item $p_q^*(\E^{\plin|_{\underline{\mathfrak{g}}^*|_{q}}}\lambda_q)$ has constant vanishing order $\ell$ along $\field P(\underline{\mathfrak{g}}^*\at{q})$, for any volume form $\lambda_q$ on $\underline{\mathfrak{g}}^*\at{q}$ and all $q\in N$.  
    \end{itemize}
 Here, $p_q\colon \blup(\underline{\mathfrak{g}}^*\at{q},\{0_q\})\to \underline{\mathfrak{g}}^*\at{q}$ is the blowdown map of the fibre over $q$, which can be seen as the restriction of the blowdown map $p\colon \blup(\underline{\mathfrak{g}}^*,N)\to \underline{\mathfrak{g}}^*$ to $p^{-1}(\underline{\mathfrak{g}}^*\at{q})$. 
\end{lemma}

\begin{proof}
The statements are local on $N$ and do not depend on the chosen volume forms. So we may assume that $\underline{\mathfrak{g}}^*=\mathbb{R}^m\times N$, and choose a product volume form $\lambda= \lambda_{\mathbb{R}^m}\wedge \lambda_N$. For any $q\in N$ we have that ${\plin}\at{\mathbb{R}^m\times\{q\}}$ is a linear Poisson structure ${\pi}_q$ on $\mathbb{R}^m$. With this, we have that
\begin{equation}\label{eq:puling_to_fibre}
p^*(\E^{\plin}\lambda)\at{\mathbb{R}^m\times \{q\}}=p_q^*(\E^{{\pi}_q}\lambda_{\mathbb{R}^m})\wedge\lambda_{N}.
\end{equation}

If $p^*(\E^{\plin}\lambda)$ has constant vanishing order $\ell$ along $\field P(\underline{\mathfrak{g}}^*)=\field P(\mathbb{R}^m)\times N$, then \[\tilde{x}_i^{-\ell}p^*(\E^{\plin}\lambda)\at{U_{i}\times N}\] is smooth and nowhere vanishing, where $(U_i,(\tilde{x}_1,\ldots,\tilde{x}_m))$ is a standard chart on $\blup(\mathbb{R}^m,\{0\})$ with
$\field P(\mathbb{R}^m)\cap U_{i}=\{\tilde{x}_i=0\}$. 
Fix $q\in N$. Restricting to $\mathbb{R}^m\times \{q\}$, \eqref{eq:puling_to_fibre} implies that 
\[\tilde{x}_i^{-\ell}p_q^*(\E^{{\pi}_q}\lambda_{\mathbb{R}^m})\at{U_i}\] is smooth and nowhere vanishing. Hence, $p_q^*(\E^{{\pi}_q}\lambda_{\mathbb{R}^m})$ has constant vanishing order $\ell$ along $\field P(\mathbb{R}^m)$. 

The converse is proven exactly in the same way: if for all $q\in N$,  
$\tilde{x}_i^{-\ell}p_q^*(\E^{{\pi}_q}\lambda_{\mathbb{R}^m})\at{U_i}$ is smooth and nowhere vanishing, then, by \eqref{eq:puling_to_fibre},  also $\tilde{x}_i^{-\ell}p^*(\E^{\plin}\lambda)\at{U_{i}\times N}$ is smooth and nowhere vanishing.
\end{proof}

Next, for the dual of a Lie algebra, we relate the height to the vanishing order.

\begin{lemma}\label{lemma:lifting_submanifold_height}
    Let $\liealg g$ be a Lie algebra and denote the linear Poisson structure on $\liealg g^\ast$ by $\pi_\lin$. Fix a volume form $\lambda$ on $\liealg g^\ast$. For $\xi\in \liealg g^\ast\setminus \{0\}$, consider the one-dimensional vector space
    \[\Sigma_{[\xi]}=(\mathbb{R}\backslash\{0\})\cdot \xi\, \cup\, \{[\xi]\}\subseteq \blup(\liealg g^\ast,\{0\}).\]
The vanishing order of $p^\ast (\E^{\plin}\lambda)\at{\Sigma_{[\xi]}}$ at $[\xi]\in \field P(\liealg g^\ast)$ is $\dim \liealg g-1-\mathrm{height}(\xi)$.
\end{lemma}

Regarding the blowup as the tautological line bundle over the projective space:
\[\sigma\colon  \blup(\mathfrak{g}^*,\{0\})\to \field P (\mathfrak{g}^*),\]
the line $\Sigma_{[\xi]}=\sigma^{-1}([\xi])$ is just the fibre over $[\xi]$.

\begin{proof}[of Lemma \ref{lemma:lifting_submanifold_height}]
Fix $\xi \in \liealg g^\ast\setminus \{0\}$. Let $b_1,\dots,b_{d}$ be a basis of $\liealg g$ such that $\xi(b_1)=1$ and $\xi(b_i)=0$, for $2\leq i\leq d$. Let $\pi_{ij}^k$ denote the corresponding structure constants of $\liealg g$, i.e.\  
    \[ [b_i,b_j]=\sum_k \pi^k_{ij}\, b_k.\] 
    The height of $\xi$ is half the rank of the skew-symmetric matrix 
    \begin{equation*}
        M_1:=(\pi_{ij}^1)_{i,j\neq 1}=(\xi[b_i,b_j])_{i,j\neq 1}.
    \end{equation*}
     Indeed, the height of $\xi$ is precisely half the rank of the skew-symmetric bilinear form $\DEC \xi$ restricted to ${\ker \xi}$. In the basis $\{b_i\}_{i=2}^d$ this form is represented by the matrix $(\DEC\xi)(b_i,b_j)=-\xi([b_i,b_j])$.

Let $\{x_i\}_{i=1}^d$ denote the linear coordinates on $\mathfrak{g}^*$ induced by the basis $\{b_i\}_{i=1}^d$. In these coordinates, the Poisson structure is given by
    \begin{equation*}
        \pi_\lin =  \frac{1}{2}\sum_{i,j,k} x_k \pi^k_{ij}\frac{\partial}{\partial x_i}\wedge\frac{\partial}{\partial x_j}.
    \end{equation*}
    Consider the chart $(U_1,\tilde{x})$ from \eqref{eq:blowdowncoord} on $\blup(\mathfrak{g}^*,\{0\})$. Using that $p^\ast x_1=\tilde{x}_1$, $p^\ast x_j=\tilde{x}_1 \tilde{x}_j$ for $j\neq 1$, and \eqref{eq:vfliftcoords}, we see that the coefficient functions of $\tilde{\pi}=p^\ast \pi_\lin \at{U_1\setminus \field P(\mathfrak{g}^*)}$ are given by
    \begin{equation*}
        \frac{1}{2}\tilde{\pi}_{ij}=\begin{cases}
            \tfrac{1}{2\tilde{x}_1} \big[ \pi_{ij}^1-\tilde{x}_j\pi_{i1}^1-\tilde{x}_i \pi_{1j}^1 +\sum_{k\neq 1} \tilde{x}_k \big(  \pi_{ij}^k-\tilde{x}_j\pi_{i1}^k-\tilde{x}_i \pi_{1j}^k   \big) \big]& \text{if } i,j\neq 1\\
             \tfrac{1}{2}(\pi_{ij}^1+\sum_{k\neq 1}\tilde{x}_k \pi_{ij}^k) &\text{if } i=1 \text{ or }j=1.
        \end{cases}
    \end{equation*}
Note that the line $\Sigma_{[\xi]}$ is included in $U_1$, and it is given by $\Sigma_{[\xi]}=\{\tilde{x}_2=\ldots=\tilde{x}_{d}=0\}$. Along this line, the coefficient functions simplify to
    \begin{equation*}
        \tilde{\pi}_{ij}=\begin{cases}
            \tfrac{1}{\tilde{x}_1} \pi_{ij}^1& \text{if } i,j\neq 1\\
             \pi_{ij}^1&\text{if } i=1 \text{ or }j=1.
        \end{cases}
    \end{equation*}
Then we can write 
\begin{equation}\label{eq:decompose_pi_above}
\tilde{\pi}=\frac{1}{\tilde{x}_1}\tilde{M}_1+\frac{\partial}{\partial \tilde{x}_1}\wedge v, \quad \textrm{where}\ v=\sum_{j=2}^d\pi_{1,j}^1\frac{\partial}{\partial \tilde{x}_j},
\end{equation}
and $\tilde{M}_1$ is the matrix $M_1$ viewed as a constant bivector field. 

   Let $\lambda=\D x_1\wedge \ldots \wedge \D x_{d}$. We calculate the restriction to $\Sigma_{[\xi]}$ of the pulled back spinor.
        \begin{align*} p^\ast(\E^{\plin}\lambda)\at{\Sigma_{[\xi]}}=&\E^{p^\ast(\pi_\lin)}p^\ast\lambda\at{\Sigma_{[\xi]}}\\
         =&\tilde{x}_1^{d-1}\E^{\tilde{\pi}}\D \tilde{x}_1\wedge\cdots\wedge \D\tilde{x}_{d}\\
         =& \sum_{k=0}^{\lfloor \tfrac{d}{2}\rfloor} \frac{1}{k!}\tilde{x}_1^{d-1-k}(\I_{\tilde{x}_1 \tilde{\pi}})^k\D \tilde{x}_1\wedge\cdots\wedge \D\tilde{x}_{d}\\
         =& \sum_{k=0}^{\lfloor \tfrac{d}{2}\rfloor} \frac{1}{k!}\tilde{x}_1^{d-1-k}(\I_v+\D \tilde{x}_1\wedge)(\I_{\tilde{M}_1})^k\D \tilde{x}_2\wedge\cdots\wedge \D\tilde{x}_{d},
        \end{align*}
where in the last step, we used \eqref{eq:decompose_pi_above}. In this sum, each power $\tilde{x}_1^{d-1-k}$ is multiplied with a constant form, which is non-zero precisely when $(\I_{\tilde{M}_1})^k\D \tilde{x}_2\wedge\cdots\wedge \D\tilde{x}_{d}\neq 0$, which is equivalent to $\mathrm{rank}(M_1)\geq 2k$. Because $\tfrac{1}{2}\mathrm{rank}(M_1)=\mathrm{height}(\xi)$, the largest $k$ such that the coefficient of $\tilde{x}_1^{d -1-k}$ in $p^\ast(\E^{\plin}\lambda)\at{\Sigma_{[\xi]}}$ is non-zero is $\mathrm{height}(\xi)$. Hence, the vanishing order is 
$\ell=d-1-\mathrm{height}(\xi)$.
\end{proof}

The next lemma will conclude the proof of Theorem \ref{theorem:summary_submanifold}. 

\begin{lemma}
    Let $\liealg g$ be a Lie algebra and denote the linear Poisson structure on $\liealg g^\ast$ by $\pi_\lin$. Fix a volume form $\lambda$ on $\liealg g^\ast$. The following are equivalent. 
    \begin{itemize}
    \item The pullback form $p^*(\E^{\plin}\lambda)$ has constant vanishing order $\ell$ along $\field P(\mathfrak{g}^*)$.
    \item The height of elements $\xi\in \mathfrak{g}^*\backslash \{0\}$ is a constant $k$.
    \end{itemize}
    In this case, $\ell+k=\dim\mathfrak{g}-1$.
\end{lemma}
\begin{proof}
In view of the previous lemma, we need to check that the following are equivalent:
\begin{itemize}
\item $p^*(\E^{\plin}\lambda)$ has constant vanishing order $\ell$ along $\field P (\mathfrak{g}^*)$, 
\item for each $[\xi]\in \field P (\mathfrak{g}^*)$, $p^*(\E^{\plin}\lambda)\at{\Sigma_{[\xi]}}$ has vanishing order $\ell$ at $[\xi]$.
\end{itemize}
Let $(U_i,\tilde{x})$ be a standard coordinate system on $\blup(\mathfrak{g}^*,\{0\})$, with $U_i\cap \field P (\mathfrak{g}^*)=\{\tilde{x}_i=0\}$. The first condition is equivalent to $\tilde{x}_i^{-\ell}p^*(\E^{\plin}\lambda)$ being smooth and nowhere vanishing on $U_i$. 
This implies that the restriction $\tilde{x}_i^{-\ell}p^*(\E^{\plin}\lambda)\at{\Sigma_{[\xi]}}$
is smooth and nowhere vanishing. For any $[\xi]\in U_i$, we have that $\tilde{x}_i\colon \Sigma_{[\xi]}\to \RR$ is a linear isomorphism. Hence, $p^*(\E^{\plin}\lambda)\at{\Sigma_{[\xi]}}$ has vanishing order $\ell$ at $[\xi]$. The converse implication is proven in exactly the same way. 
\end{proof}

\section{Lifting Poisson structures: a geometric approach}\label{ref:sec:lift_Poisson}

In this section, we present a more geometric proof, which does not use spinors, of Theorem \ref{theorem:lifting_by_spinors} for the particular case of blowing up a zero of a Poisson manifold. This implies also the case of symplectic leaves of a Poisson manifold and, more generally, of presymplectic leaves of a Dirac manifold, by using the splitting theorem of Weinstein \cite[Theorem 2.1]{weinstein83a} and Blohmann \cite{blohmann:2017a}, respectively. Moreover, we give alternative descriptions of the constant height condition for a Lie algebra. 

\begin{theorem}\label{theorem:point_lift}
    Let $(M,\pi)$ be a Poisson manifold and $q\in M$ a zero of $\pi$, i.e.\ $\pi(q)=0$. Let $\liealg g:=T^*_qM$ denote the isotropy Lie algebra of $\pi$ at $q$.
    The following are equivalent.
    \begin{enumerate}
        \item \label{theorem:point_lift:item:lifts} The Dirac structure $\graph(\pi)$ lifts to a Dirac structure $\tilde
        L$ on $\blup(M,\{q\})$.
        \item \label{theorem:point_lift:item:distribution} There exists $ k\in \field N_0 $ such that for all $v\in T_qM\setminus \{0\}$ the subspace
        \begin{equation}\label{eq:assvcirc}
            D_{[v]}:=\big\{ (\widetilde{ \pi^\sharp \alpha })_{[v]}\colon \alpha\in \Omega^1(M),\  \alpha_q(v)=0 \big\}\subseteq T_{[v]} \field P(T_qM)
        \end{equation}
        has rank $2k$. Here, $\tilde{X}$ denotes the lift to $\blup(M,\{q\})$ of a vector field $X$ on $M$ vanishing at $q$. 
        \item \label{theorem:point_lift:item:coadjoint}  There exists $ k\in \field N_0 $ such that for all $\xi\in \liealg g^\ast\setminus\{0\}$,  the coadjoint orbit $\mathcal{O}_\xi$  through $\xi$ satisfies
        \begin{equation}
            \dim(\mathcal{O}_\xi)=\begin{cases}
            2k+2 &\text{ if }T_\xi\mathcal{O}_\xi\text{ contains the radial line }\field R \xi,\\
            2k &\text{ otherwise}.
            \end{cases}
        \end{equation}
    \item \label{theorem:point_lift:item:height} The Lie algebra $\g$ has constant height $k$, i.e.\ for all $\xi\in \liealg g^\ast\setminus\{0\}$ we have
    \begin{equation}
        \xi\wedge (\DEC\xi)^{k}\neq 0 \quad \text{ and }\quad \xi\wedge (\DEC\xi)^{k+1}=0,
    \end{equation}
     where $\DEC\colon \wedge^\bullet \liealg g^\ast\to \wedge^{\bullet+1} \liealg g^\ast$ denotes the Chevalley-Eilenberg differential.
    \end{enumerate}
If any (and thus all) conditions are satisfied, the blowdown map is a forward Dirac map, and
\begin{equation*}
    \tilde{L}\cap T\field P(T_qM) = D. 
\end{equation*} 
Moreover, $\tilde{L}$ is the graph of a Poisson structure if and only if $k=0$.
\end{theorem}

\begin{remark}
The integers $k$ appearing in conditions 2,3 and 4 all agree.
\end{remark}

\begin{remark}\label{rem:D_using_lin}
The space $D$ appearing in condition 2 depends only on the linearisation $\plin$ of $\pi$ at $q$. Indeed, in a chart centred at $q$, the vector fields $(\pi-\plin)^{\sharp}\alpha$ vanish to second order at $0$, and so their lift to $\blup(M,\{q\})$ vanishes along $\mathbb{P}(T_qM)$.
\end{remark}

\subsection{Examples for Theorem \ref{theorem:point_lift}}
 
Before proving Theorem \ref{theorem:point_lift}, we present some explicit examples.

\begin{example}\label{ex:so3}
We check explicitly the conditions of Theorem \ref{theorem:point_lift} for the Poisson manifold $\liealg{so}(3)^*$, with the linear Poisson structure $\pi$, showing that it lifts to a Dirac structure on the blowup which is not a Poisson structure, and describing its geometry. We start from condition 4.

\begin{itemize}
    \item[4.] 
The Lie algebra $\liealg{so}(3)$ admits a basis $\{X_1,X_2,X_3\}$ such that $[X_1,X_2]=X_3$, $[X_2,X_3]=X_1$, $[X_3,X_1]=X_2$. Denote by  $\{\theta_1,\theta_2,\theta_3\}$ the dual basis.
For all
$\xi\in \liealg{so}(3)^*\setminus\{0\}$ we have
$$\xi\wedge (\D_{\liealg{so}(3)}\xi)=-\big(\sum_{j=1}^{3} (\xi(X_j))^2\big) \theta_1\wedge \theta_2\wedge \theta_3\neq 0,$$ while  
$\xi\wedge (\D_{\liealg{so}(3)}\xi)^2=0$ for dimensional reasons. Hence, condition 4 is satisfied for $k=1$.

\item[3.] The coadjoint orbits (apart from the origin) are concentric spheres around the origin. In particular, they have dimension $2$, and they do not contain any radial line. Therefore, condition 3 is satisfied for  $k=1$.

\item[1.] Denote by $\{x_i\}$ the linear coordinates on $\liealg{so}(3)^*$ arising from the basis $\{X_i\}$. A frame for $\graph(\pi)$ is given by $\{(\pi^\sharp \D x_i, \D x_i)\}$, where for instance $\pi^\sharp \D x_1 =x_3{\partial}_2-x_2{\partial}_3$. In the chart $U_1$ of Section \ref{sec:blowup} one has that $(\widetilde{\pi^\sharp \D x_i},p^*\D x_i)$ equals, for $i=1,2,3$ respectively,
\begin{equation*}
\begin{pmatrix}\tilde{x}_3\tilde{\partial}_2-\tilde{x}_2\tilde{\partial}_3\\\D \tilde{x}_1
\end{pmatrix},
\quad\ 
\begin{pmatrix}-\tilde{x}_1\tilde{x}_3\tilde{\partial}_1+\tilde{x}_2\tilde{x}_3\tilde{\partial}_2+(1+\tilde{x}_3^2)\tilde{\partial}_3\\
\tilde{x}_1 \D \tilde{x}_2+\tilde{x}_2 \D \tilde{x}_1
\end{pmatrix},
\quad\ 
\begin{pmatrix}\tilde{x}_1\tilde{x}_2\tilde{\partial}_1-(1+\tilde{x}_2^2)\tilde{\partial}_2-\tilde{x}_2\tilde{x}_3\tilde{\partial}_3\\
\tilde{x}_1 \D \tilde{x}_3+\tilde{x}_3 \D \tilde{x}_1
\end{pmatrix}.
\end{equation*}
This can be computed using \eqref{eq:blowdowncoord} and \eqref{eq:vfliftcoords}; here we use the notation $\tilde{\partial}_j:={\partial}_{\tilde{x}_j}$. We know that these three sections are linearly independent away from $\field P (\mathfrak{so}(3)^*)\cap U_1= \{\tilde{x}_1=0\}$, but they are also on $\{\tilde{x}_1=0\}$, since there the determinant of the coefficients of $\tilde{\partial}_2, \tilde{\partial}_3, \D \tilde{x}_1$ is given by $(1+
(\tilde{x}_2)^2+(\tilde{x}_3)^2)^2>0$.
 Hence, these three sections span a Dirac structure on $U_1$ that lifts $\pi$. A similar computation can be done for the charts $U_2$ and $U_3$, showing directly that $\pi$ lifts to a Dirac structure $\tilde{L}$.
\item[2.]
Let $v=(v_1,v_2,v_3)\in \mathfrak{so}(3)^*\backslash \{0\}$. Without loss of generality, assume that $v_1\neq 0$. Then $[v]$ is covered by the chart $(U_1,\tilde{x})$, and in these coordinates $(\tilde{x}_1([v]),\tilde{x}_2([v]),\tilde{x}_3([v])) =(0,\frac{v_2}{v_1},\frac{v_3}{v_1})$. A basis for the covectors that annihilate $v$ is given by
\[\alpha_2:=v_3\D x_1-v_1\D x_3, \quad \alpha_3:=v_1\D x_2-v_2\D x_1,\]
which we regard also as constant 1-forms.
Using the above expressions for $\widetilde{\pi^{\sharp}\D x_i}$, we obtain
\[(\widetilde{\pi^{\sharp}\alpha_2})_{[v]}=\frac{v_1^2+v_2^2+v_3^2}{v_1}\tilde{\partial}_{2}, \qquad (\widetilde{\pi^{\sharp}\alpha_3})_{[v]}=\frac{v_1^2+v_2^2+v_3^2}{v_1}\tilde{\partial}_{3}.\]
Hence, we obtain that
\[D_{[v]}=\Span\{\tilde{\partial}_{2},\tilde{\partial}_{3}\}=T_{[v]}\field P (\mathfrak{so}(3)^*).\]
Hence, $D=T\field P (\mathfrak{so}(3)^*)$, and so the rank of $D$ is two  everywhere. So 2.\ holds for $k=1$.
\end{itemize}
 Further, it is easy to see that on $\field P (\mathfrak{so}(3)^*)\cap U_1= \{\tilde{x}_1=0\}$ the   Dirac structure $\tilde{L}$ is spanned also by $(\tilde{\partial}_2,0),(\tilde{\partial}_3,0), (0,   \D \tilde{x}_1)$. This implies that $\field P (\mathfrak{so}(3)^*)$ together with the zero 2-form is a presymplectic leaf of $\tilde{L}$. In particular, the Dirac structure $\tilde{L}$ is regular.
\end{example}

\begin{example}\label{ex:sl2}
We briefly look at 
  the conditions of Theorem \ref{theorem:point_lift} for the Poisson manifold $(\liealg{sl}_2(\field R)^*,\pi)$, building on Example \ref{ex:so3}, and showing that the Poisson structure  does not lift.
\begin{itemize}
    \item[4.]  
For all
$\xi\in\liealg{sl}_2(\field R)^*\setminus\{0\}$ we have
$$\xi\wedge (\D_{\liealg{sl}_2(\field R)}\xi)=\big(-(\xi(e_1))^2 -(\xi(e_2))^2 +(\xi(e_3))^2 \big)\, e_1^*\wedge e_2^*\wedge e_3^*,$$ which vanishes exactly on the two coadjoint orbits given by half-cones.
\item[3.] The coadjoint orbits (apart from the origin)  are given by hyperboloids, paraboloids and two half-cones. They all have  dimension $2$, but the half-cones contain the radial line, while the other orbits do not. 
\end{itemize}
\end{example}

\subsection{Proof of the equivalence $\ref{theorem:point_lift:item:lifts} \Leftrightarrow \ref{theorem:point_lift:item:distribution}$ in Theorem \ref{theorem:point_lift}}

We prove now the equivalence $\ref{theorem:point_lift:item:lifts} \Leftrightarrow \ref{theorem:point_lift:item:distribution}$,
together with the two final statements of Theorem \ref{theorem:point_lift}.

We may assume that $M=\mathbb{R}^n$, with coordinates $(x_1,\dots,x_n)$, and that $q=0$.  
 We use the chart $(U_n,\tilde{x})$ of $\blup(\mathbb{R}^n, \{0\})$ from Section \ref{sec:blowup}, in which $U_n\cap \field P(\mathbb{R}^n) =\{\tilde{x}_n=0\}$.

 Note that for all $i<n$ we have
    \begin{equation}\label{eq:howto_x_n_d_x_i}
        p^\ast \D x_i- \tilde{x}_i p^\ast \D x_n= \tilde{x}_n\D \tilde{x}_i+\tilde{x}_i  \D \tilde{x}_n - \tilde{x}_i \D \tilde{x}_n=\tilde{x}_n\D \tilde{x}_i.
    \end{equation}

    We consider the sections
    $$\left\{\begin{pmatrix}
\widetilde{\pi^{\sharp} \D x_i}\\
p^*\D x_i
   \end{pmatrix}\right\}_{i\leq n}$$ 
   of $\field T \blup(\mathbb{R}^n,\{0\})$, which away from $\field P(\mathbb{R}^n)$ are a frame for the Dirac structure corresponding to $\graph(\pi)$ under the blowdown map $p$. Over $U_n$, for $i<n$, replace the $i$-th section by itself minus $\tilde{x}_i$ times the $n$-th section. This yields another family of sections with the same property,
   \begin{equation}\label{eq:frame_on_blup:coordinate_edition}
   e_n=\begin{pmatrix}
		\widetilde{\pi^{\sharp} \D x_n}\\
 \D\tilde{x}_n
			\end{pmatrix}\quad
 \text{ and }\quad
			  e_i=
     \begin{pmatrix}
\widetilde{\pi^{\sharp}\D x_i}-\tilde{x}_i\widetilde{\pi^{\sharp} \D x_n}\\
\tilde{x}_n\D \tilde{x}_i
			\end{pmatrix}\quad\text{ for }i< n,
\end{equation}
where we used \eqref{eq:howto_x_n_d_x_i} for the cotangent part. Notice that the vector component is tangent to $\field P(\mathbb{R}^n)$ by Lemma \ref{lemma:lifting_vf_to_the_blowup}, since $\pi$ vanishes at the origin, and for $i<n$ the covector components vanish on $\field P(\mathbb{R}^n)$. 

\begin{lemma}\label{lem:claim}For any $[v]\in \field P(\mathbb{R}^n)\cap U_n$, the span of the tangent components of the $\{e_i\}_{i<n}$ is precisely the subspace $ D_{[v]}$ given in \eqref{eq:assvcirc}.
\end{lemma}
\begin{proof}
Indeed, fix $[v] \in U_n\cap \field P(\mathbb{R}^n)$ and set $\alpha_i^{[v]}:= \D x_i-\tilde{x}_i([v]) \D x_n$. Then
\begin{equation*}
\big(\widetilde{\pi^\sharp \D x_i}-\tilde{x}_i \widetilde{\pi^\sharp \D x_n}\big) _{[v]}=
    (\widetilde{\pi^\sharp \D x_i})_{[v]} -\tilde{x}_i([v]) (\widetilde{\pi^\sharp \D x_n})_{[v]}=(\widetilde{\pi^\sharp \alpha_i^{[v]}})_{[v]}.
\end{equation*}
The key observation is that $\alpha_i^{[v]}$ annihilates $v\in T_0 \mathbb{R}^n$, as $\tilde{x}_i([v])=\tfrac{v_i}{v_n}$. For $i<n$ these covectors form a basis of the annihilator $(\field R v)^\circ$, proving the claim.
\end{proof}

\subsubsection{Proof of the implication $\ref{theorem:point_lift:item:lifts} \Leftarrow \ref{theorem:point_lift:item:distribution}$ in Theorem \ref{theorem:point_lift}}

Assume the constant rank condition in the statement, i.e.\ the subspaces $D_{[v]}$, $[v]\in \field P(\mathbb{R}^n)$, have constant rank $2k$. 

Let $[v_0]\in\PP(\RR^n)\cap U_n$. Using Lemma \ref{lem:claim}, after restricting to a neighbourhood $U'\subseteq U_n$ of $[v_0]$, by permuting the $\{e_i\}_{1\leq i\leq  n-1}$, we may assume that
\begin{itemize}
    \item the sections $\{e_i\}_{1\leq i\leq 2k}$ restricted to $\PP(\RR^n)\cap U'$ are a frame for the distribution $D\at{\PP(\RR^n)\cap U'}$.
\end{itemize}
The other sections can be written along $\PP(\RR^n)\cap U'$ as unique linear combinations: 
\[e_j\at{\PP(\RR^n)\cap U'}=\sum_{i=1}^{2k}a_j^i e_i\at{\PP(\RR^n)\cap U'},\quad 2k+1 \leq j\leq n-1.\]
Extend the smooth functions $a_j^i$ to $U'$ (and denote these in the same way), and define the sections
\[e'_j:=e_j-\sum_{i=1}^{2k}a_j^i e_i,\quad 2k+1 \leq j\leq n-1.\]
We obtain sections of 
$\TT U'$ 
$$ e_1,\dots,e_{2k},e'_{2k+1},\dots,e'_{n-1},e_n,$$
with the following properties: 
\begin{itemize}
\item on $U'\setminus \PP(\RR^n)$ they form a frame for the lifted Dirac structure,
    \item the  $\{e_j'\}_{2k+1\leq  j\leq n-1}$ vanish on $\PP(\RR^n)\cap U'$. 
\end{itemize}
Since $\tilde{x}_n$ is a defining function for $\field P(\mathbb{R}^n)$ inside $U_n$, we can divide the  $e_j'$ by $\tilde{x}_n$ to obtain a smooth section of $\field T U'$, and obtain sections
\begin{equation}\label{eq:frametildeL:coordinate_edition}
    e_1,\dots,e_{2k},\frac{1}{\tilde{x}_n}e'_{2k+1},\cdots,\frac{1}{\tilde{x}_n}e'_{n-1},e_n
    \end{equation}
which on $U'\backslash \PP(\RR^n)$ are a frame for the lifted Dirac structure. These sections \eqref{eq:frametildeL:coordinate_edition} are pointwise linearly independent also at points of $\field P(\mathbb{R}^n)\cap U'$. Indeed, using the explicit form \eqref{eq:frame_on_blup:coordinate_edition}, along $\field P(\mathbb{R}^n)\cap U'$ these elements have the following properties:
\begin{itemize}
    \item $e_i\at{\field P(\mathbb{R}^n)\cap U'}$, for $1\leq i \leq 2k$, have only vector components and are linearly independent;
    \item $\tfrac{1}{\tilde{x}_n}e'_{j}\at{\field P(\mathbb{R}^n)\cap U'}$, $2k+ 1\leq j\leq n-1$, has cotangent components
    \[\Big(\D\tilde{x}_j-\sum_{i=1}^{2k}a_j^i \D\tilde{x}_i\Big)\at{\field P(\mathbb{R}^n)\cap U'};\]
    \item $e_n$ has cotangent component $\D \tilde{x}_n\at{\field P(\mathbb{R}^n)\cap U'}.$
\end{itemize}

By covering $\PP(\RR^n)$ with such open sets $U'$, we obtain that the vector subbundle \[p^\ast(\graph(\pi)\at{\mathbb{R}^n\setminus \{0\}})\subseteq \field T\blup(\mathbb{R}^n, \{0\})\at{\blup(\mathbb{R}^n,\{0\})\setminus \field P(\mathbb{R}^n)}\] extends to a smooth vector subbundle of $\field T \blup(\mathbb{R}^n,\{0\})$. By \cite[Remark 2.9.]{blohmann:2017a}, the extension is a Dirac structure, as all axioms hold on the dense open subset $\blup(\mathbb{R}^n,\{0\})\setminus \field P(\mathbb{R}^n)$. 

It is a good moment to prove one of the additional, final statements of Theorem \ref{theorem:point_lift}. 

\begin{lemma}
The pullback Dirac structure is Poisson if and only if $k=0$.
\end{lemma}
\begin{proof}
For the frame \eqref{eq:frametildeL:coordinate_edition} of the lift $\tilde{L}$, the cotangent components of the first $2k$-elements vanish along $\PP(\RR^n)$ and the cotangent components of the last $n-2k$ are linearly independent along $\PP(\RR^n)$. This proves the claim. 
\end{proof}

\subsubsection{Proof of the implication $\ref{theorem:point_lift:item:lifts} \Rightarrow \ref{theorem:point_lift:item:distribution}$ in Theorem \ref{theorem:point_lift}}

Assume that $ \graph(\pi) $ lifts to a Dirac structure $\tilde{L}$ on $ \blup(\RR^m,\{0\})$. 

The crucial observation for the proof of the implication $\ref{theorem:point_lift:item:lifts} \Rightarrow \ref{theorem:point_lift:item:distribution}$
is the following equality, asserted in Theorem \ref{theorem:point_lift}: 
\begin{equation}\label{eq:tLD}
     \tilde{L}\at{\PP(\RR^n)}\cap T \PP(\RR^n)=D.
\end{equation}
Note that the rank of the left-hand side (the intersection of two subbundles) can locally only decrease, and the rank of the right-hand side (the image of a vector bundle map) can locally only increase. Therefore, their equality implies \ref{theorem:point_lift:item:distribution} in Theorem \ref{theorem:point_lift}: $D$ has constant rank.

To prove \eqref{eq:tLD}, we start with a technical result.

\begin{lemma}\label{lem:alpha0} Let $\gamma\in \Omega^1(\blup(\RR^n,\{0\}))$. If $\gamma$ vanishes at $[v]\in \PP(\RR^n)$, then there exists $\beta\in \Omega^1(\RR^n)$ such that $p^*(\beta_{tv})=\gamma_{p^{-1}(tv)}$ for all $t\neq 0$. Moreover, $\beta_0(v)=0$.
\end{lemma}

\begin{proof} Without loss of generality, we may assume that $v_n\neq 0$, where $v=(v_1,\ldots,v_n)$. We will use coordinates $\tilde{x}_1,\dots,\tilde{x}_n$ on $U_n\subseteq \blup(\RR^n,\{0\})$ as in Section \ref{sec:blowup}. The blowdown map $p$ satisfies the following, as one sees using \eqref{eq:blowdowncoord}: 
\begin{equation*}
    p^\ast (\D x_n)=\D \tilde{x}_n \quad\text{and}\quad p^\ast\Big( \frac{\D x_i}{x_n}-\frac{x_i}{x_n^2}\D x_n\Big) =\D \tilde{x}_i  \text{ for all } 1\leq i\leq n-1.
\end{equation*}
On the domain $U_n$ of the chart,  write $\gamma$ as $\sum_{i=1}^n f_i \D\tilde{x}_i$ for smooth functions $f_i$. Consider the line
\begin{equation}\label{line:sigma_v}
    \Sigma_{[v]}:=p^{-1}(\field Rv\setminus \{0\})\cup [v]\subseteq \blup(\RR^n,\{0\}).
\end{equation}
The blowdown map gives an isomorphism $\Sigma_{[v]}\cong \mathbb{R}v$. 
Since $\gamma$ vanishes at $[v]$, on $\Sigma_{[v]}$ we can write $f_i=\tilde{x}_n\widehat{f}_i\circ p$ for all $i$, for smooth functions $\widehat{f}_i$ on $\RR\cdot v$. We extend these functions to smooth functions on $\RR^n$, which we denote in the same way. 
Using \eqref{eq:blowdowncoord}, for $t\neq 0$ we have 
\begin{equation*}
\begin{aligned}
\gamma\at{p^{-1}(tv)}&=p^*\Big(\sum_{i<n}x_n \widehat{f}_i(x)\Big(\frac{\D x_i}{x_n}-\frac{x_i}{x_n^2}\D x_n\Big)+x_n \widehat{f}_n(x) \D x_n \Big)\Big|_{p^{-1}(tv)}\\
&=p^*\Bigg(\sum_{i< n} \widehat{f}_i(tv) \D x_i\at{tv}+
\Big(tv_n \widehat{f}_n(tv) - \sum_{i< n} \widehat{f}_i(tv) \frac{v_i}{v_n}\Big)\D x_n\at{tv}
\Bigg).
\end{aligned}
\end{equation*}
This expression extends smoothly to $t=0$, and yields an explicit formula:
\[\beta:=\sum_{i< n} \widehat{f}_i(x) \D x_i+
\Big(x_n \widehat{f}_n(x) - \sum_{i< n} \widehat{f}_i(x) \frac{v_i}{v_n}\Big)\D x_n.\]
At $t=0$, we have that:
\[\beta_0(v)=\sum_{i< n} \widehat{f}_i(0) v_i+ \Big(0
- \sum_{i< n} \widehat{f}_i(0) \frac{v_i}{v_n}\Big)v_n=0.\]
\end{proof}

\begin{proof}[of Equation \eqref{eq:tLD}]
We just need to prove the inclusion ``$\subseteq $'' in \eqref{eq:tLD}, since the other inclusion follows from Lemma \ref{lem:claim} and 
the lines preceding it.
Let $[v]\in \PP(\RR^n)$ and let $X_{[v]}\in \tilde{L}_{[v]} \cap T_{[v]}\PP(\RR^n)$. Extend $(X_{[v]},0)$ to a section $(X,\gamma)\in \Gamma^\infty(\tilde{L})$. 
Then there exists a unique $\alpha\in \Omega^1(\RR^n \setminus\{0\})$ such that $$\begin{pmatrix} X\\
\gamma \end{pmatrix}\Big|_{\blup(\RR^n,\{0\})\setminus\PP(\RR^n)}=\begin{pmatrix} \widetilde{\pi^{\sharp} \alpha}\\
p^*\alpha
\end{pmatrix}
\Big|_{\blup(\RR^n,\{0\})\setminus\PP(\RR^n)},$$ since the blowdown map $p$ is a diffeomorphism away from $\PP(\RR^n)$. Lemma \ref{lem:alpha0} implies that there exists $\beta\in \Omega^1(\RR^n)$ which agrees with $\alpha$ on the line $\RR v\setminus\{0\}$. We consider a limit along the line $\Sigma_{[v]}$ \eqref{line:sigma_v}. 
We have 
$$X_{[v]}=\lim_{t\to 0} \widetilde{\pi^{\sharp} \alpha}\at{p^{-1}(tv)}=\lim_{t\to 0} \widetilde{\pi^{\sharp} \beta}\at{p^{-1}(tv)}= \widetilde{\pi^{\sharp} \beta}\at{[v]}=\widetilde{\pi^{\sharp} \beta_0}\at{[v]},$$
where $\beta_0$ is viewed as a constant 1-form, and in the last equality uses that $\pi^{\sharp}(\beta-\beta_0)$ vanishes to second order at $0$, so its lift is zero along $\mathbb{P}(\RR^n)$. By Lemma \ref{lem:alpha0}, $\beta_0(v)=0$.
Hence, $X_{[v]}\in D_{[v]}$.
\end{proof}

Finally, we prove also the following assertion of Theorem \ref{theorem:point_lift}. 

\begin{lemma}The blowdown map $\blowdown{p}{\RR^n}{\{0\}}$ is a forward Dirac map. 
\end{lemma}
\begin{proof}
We need to show that $p_{!}(\tilde{L}_{[v]})=T^*_0\RR^n$, for all $[v]\in \PP(\RR^n)$. Let $\alpha_0\in T^*_0\RR^n$. We regard $\alpha_0$ as a constant 1-form. Then $(\widetilde{\pi^{\sharp}\alpha_0},p^*\alpha_0)$ belongs to $\tilde{L}$ outside of $\PP(\RR^n)$. Since $\tilde{L}$ is a closed subset, also $(\widetilde{\pi^{\sharp}\alpha_0},p^*\alpha_0)_{[v]}\in \tilde{L}_{[v]}$, for all $[v]\in \PP(\RR^n)$. Since $T_{[v]}p((\widetilde{\pi^{\sharp}\alpha_0})_{[v]})=0$, $\alpha_0\in p_{!}(\tilde{L}_{[v]})$.
\end{proof}

\subsection{The equivalence $\ref{theorem:point_lift:item:distribution}\Leftrightarrow\ref{theorem:point_lift:item:coadjoint}$ in Theorem \ref{theorem:point_lift}}

We now give a more explicit description of the singular distribution $D$ in  \eqref{eq:assvcirc}, solely in terms of the coadjoint orbits of the Lie algebra $\g$.
\begin{lemma}\label{lemma:Distribution_and_coadj_orbits}
Let $\g$ be a Lie algebra. Let $\xi\in \g^*$ be a non-zero element, and denote by $\cO_\xi$ the coadjoint orbit through $\xi$. Denote by $\pi_{\lin}$ the linear Poisson structure on $\g^*$, hence $\image((\pi_{\lin})^{\sharp}_\xi)=T_\xi\cO_\xi$. Then
\begin{itemize}
    \item [i)] $T_\xi\cO_\xi$ contains the radial line $\RR \xi$ $\Leftrightarrow$ $(\pi_{\lin})^{\sharp}_\xi((\RR \xi)^{\circ})\subsetneq T_\xi\cO_\xi$, 
\end{itemize}
and in this case, we have that
\begin{itemize}
    \item [ii)] $\RR \xi\subseteq (\pi_{\lin})^{\sharp}_\xi((\RR \xi)^{\circ})$.
\end{itemize}
\end{lemma}
\begin{proof}
Since $B:=({\plin})^{\sharp}_{\xi}$ is skew-symmetric, $(\image B)^{\circ}=\ker B$. So, we obtain i):
\[
 \xi\in  \image B  \, \Leftrightarrow\,  
\RR \xi\cap \image B \neq \{0\}\, \Leftrightarrow\, 
 (\RR \xi)^{\circ} + \ker B\neq \mathfrak{g}\, \Leftrightarrow\, 
 B( (\RR \xi)^{\circ}) \neq \image B,
\]
where the last equivalence follows from the rank-nullity theorem. In this case, there is $X\in \g$ such that $B(X) =\xi$. By skew-symmetry, we get 
$0=\langle B X, X \rangle=\langle  \xi, X\rangle,$
i.e.\ $X \in (\RR \xi)^{\circ}$.
\end{proof}

From Lemma \ref{lemma:Distribution_and_coadj_orbits} we   obtain the equivalence $\ref{theorem:point_lift:item:distribution}\Leftrightarrow\ref{theorem:point_lift:item:coadjoint}$ in Theorem \ref{theorem:point_lift}.

\begin{corollary}
    Let $ (M,\pi) $ be Poisson and $ q\in M $ with $ \pi(q)=0 $, consider the Lie algebra $\g=T_q^*M$.
    For every non-zero $\xi\in \g^*$, denote by $\cO_\xi$ the coadjoint orbit through $\xi$. Then the dimension of the subspace $D_{[\xi]}$ of $T_{[\xi]}(\PP(T_{q}M))$ of  \eqref{eq:assvcirc} reads:
   $$\dim(D_{[\xi]})=
    \begin{cases}
      \dim(\cO_\xi)-2& \quad \text{if   $T_\xi\cO_\xi$ contains the radial line $\RR \xi$},\\
     \dim(\cO_\xi)& \quad \text{otherwise}.
    \end{cases}
    $$
\end{corollary}
\begin{proof}
The subspace
$\{(\pi_{\lin})^{\sharp}_\xi \alpha\colon\alpha \in (\field R \xi)^{\circ}\}$ of $\image((\pi_{\lin})^{\sharp}_\xi)=T_\xi\cO_\xi$ has codimension either zero or one. It has codimension one exactly when $T_\xi\cO_\xi$ contains the radial line $\RR \xi$, by Lemma \ref{lemma:Distribution_and_coadj_orbits} i). In that case, the image of the subspace under the projectivisation has one dimension less, by Lemma \ref{lemma:Distribution_and_coadj_orbits} ii). The statement follows by re-expressing  $D_{[\xi]}$ \eqref{eq:assvcirc} by means of $\plin$ (Remark \ref{rem:D_using_lin}).
\end{proof}

\subsection{The equivalence $\ref{theorem:point_lift:item:coadjoint}\Leftrightarrow \ref{theorem:point_lift:item:height}$ in Theorem \ref{theorem:point_lift}}    
 
The following lemma concludes the proof of Theorem \ref{theorem:point_lift}.
\begin{lemma}\label{lemma:equivalence_orbit_dga}
    Let $\g$ be a Lie algebra and $\xi\in \mathfrak{g}^*\backslash\{0\}$. Let $\mathcal{O}_{\xi}$ be the coadjoint orbit of $\xi$, and $k\in \field N_0$ be the height of $\xi$, i.e.\ 
     $\xi\wedge (\DEC\xi)^{k}\neq 0 $ and $ \xi\wedge (\DEC\xi)^{k+1}=0$. 
One of the following alternatives holds.
    \begin{enumerate}[(1)]
\item\label{type_1} $(\DEC\xi)^{k+1}=0$, which is equivalent to (1') 
$\dim(\mathcal{O}_\xi)=2k$, and also to (1'') $\xi\notin T_\xi \mathcal{O}_\xi$.
        \item\label{type_2} $(\DEC\xi)^{k+1}\neq 0$, which is equivalent to (2') $\dim(\mathcal{O}_\xi)=2k+2$, and also to (2'') $\xi\in T_\xi \mathcal{O}_\xi$.
    \end{enumerate}
    \end{lemma}

\begin{proof} 
First, note that we have 
    \begin{equation*}
        T_\xi \mathcal{O}_\xi = \big(\ker (\DEC\xi)\big)^\circ,
    \end{equation*}
    hence, 
    $$\dim (\mathcal{O}_\xi)=\rank(\DEC\xi).$$
    
\begin{enumerate}[(1)]
\item If $(\DEC\xi)^{k+1}=0$, then, since $(\DEC\xi)^{k}\neq0$, it follows that $\dim(\mathcal{O}_\xi) = \rank(\DEC\xi) =2k$. Moreover, $ \xi\notin T_\xi\mathcal{O}_\xi$. Indeed, if $\xi\in \big(\ker (\DEC\xi)\big)^\circ $, then, since $k=\tfrac{1}{2}\rank \DEC\xi$, we would obtain $\xi \wedge (\DEC\xi)^{k}=0$.
\item
    If $(\DEC\xi)^{k+1}\neq 0$, then, since 
    \begin{equation*}
        0=\DEC(0)=\DEC(\xi\wedge (\DEC\xi)^{k+1})=(\DEC\xi)^{k+2},
    \end{equation*}
    we have $\dim (\mathcal{O}_\xi)= 2k+2$. By contracting the identity $\xi\wedge (\DEC\xi)^{k+1}=0$ with any vector in  $\ker(\DEC\xi)$, we see that $\ker(\DEC\xi)\subseteq \xi^{\circ}$. Taking annihilators we obtain  $ \xi\in T_\xi\mathcal{O}_\xi$.
\end{enumerate}
Since \emph{(1')} and \emph{(2')} are mutually exclusive, it follows that they imply \emph{(1)} and \emph{(2)}, respectively. The same holds for \emph{(1'')} and \emph{(2'').} 
\end{proof}

\section{Classification of Lie algebras of constant height}\label{sec:classif}

In this section, we classify all Lie algebras that satisfy the conditions stated in Theorem \ref{theorem:lifting_by_spinors}.

\begin{theorem}\label{theorem:classification_LA_of_constant_height}
  Any Lie algebra $\liealg g$ of constant height is isomorphic to one of the following.
\begin{itemize}        \item An abelian Lie algebra $\field R^n$---this has height $0$.
 \item The semi-direct product $\field R\ltimes \field R^n$, for the representation $\lambda\mapsto \lambda\id_{\RR^n}$---this has height $0$.
        \item The Lie algebra $\liealg{so}(3)$---this has height $1$.
    \end{itemize}
\end{theorem}

Is useful to introduce the following terminology.

 \begin{definition}
     We say that an element $\xi\in \liealg g^\ast\setminus\{0\}$ is of \emph{type (1)} if it satisfies condition \ref{type_1} of Lemma \ref{lemma:equivalence_orbit_dga}, and of \emph{type (2)} if it satisfies condition \ref{type_2} of Lemma \ref{lemma:equivalence_orbit_dga}. \end{definition}

\begin{remark}\label{rem:gozeremm}
The notions of \emph{height} and \emph{type} of an element $\xi\in \liealg g^\ast\setminus \{0\}$ are encoded by the notion of \emph{Cartan class} used in 
\cite[\S 2.1]{goze.remm:2019a}. This number is characterised by (see 
\cite[\S 2.1]{goze.remm:2019a})	
\begin{align*}
			\mathrm{class}(\xi)&=2k+1 \Leftrightarrow \xi\wedge (\DEC\xi)^{k}\neq 0 \ \textrm{and}\  (\DEC\xi)^{k+1}=0
,\\
			\mathrm{class}(\xi)&=2k+2\Leftrightarrow (\DEC\xi)^{k+1}\neq 0\ \text{and}\ \xi\wedge (\DEC\xi)^{k+1}=0.\end{align*}
The definitions and $\DEC(\alpha\wedge (\DEC\alpha)^{k})=(\DEC\alpha)^{k+1}$ imply the equality:
\[\text{class}(\xi)= 2\cdot \text{height}(\xi)+\text{type}(\xi).\]
The notion of Cartan class was used in \cite[Prop. 2.9]{goze.remm:2019a} to show that $\liealg{so}(3)$ or $\liealg{sl}_2(\field R)$ are, up to isomorphism, the only Lie algebras whose non-trivial coadjoint orbits have \emph{codimension $1$}. Notice that Theorem \ref{theorem:classification_LA_of_constant_height} implies a variation of this statement: any  \emph{compact} Lie algebra, such that all coadjoint orbits (except for the origin) have the same dimension, must be abelian or isomorphic to $\liealg{so}(3)$.
\end{remark}

An important step in the proof is to show that any semisimple Lie algebra of constant height is isomorphic to $\liealg{so}(3)$ (Theorem \ref{theorem:ss_constant_height_is_so3}). In order to show compactness, we will need the following result. 

\begin{lemma}\label{lemma:ss_constant_height_is_type_1}
 For a semisimple Lie algebra $\mathfrak{g}$ of constant height, all elements in $\mathfrak{g}^*\backslash\{0\}$ are of type \ref{type_1}.
\end{lemma}
\begin{proof} Let $k$ denote the constant height of $\mathfrak{g}$. Suppose $\xi\in \liealg g^\ast\setminus\{0\}$ is of type \ref{type_2}, i.e.\ $(\DEC \xi)^{k+1}\neq 0$. Consider the complexification $\liealg g_{\field C}=\liealg g\tensor \field C$. Using $\field C$-linear extension, we regard 
 $\liealg g^\ast\subseteq \liealg g_{\field C}^\ast$. As such, our specific $\xi$ still satisfies the condition $(\DECC \xi)^{k+1}\neq 0$.    Thus, the set 
    \begin{equation*}
        U:=\{ \eta \in \liealg g^\ast_{\field C}\colon (\DECC \eta)^{k+1}\neq 0 \},
    \end{equation*}
    is non-empty, hence open and dense (because it is Zarisky open). 

    {\sf Claim:} {\it For all $\eta\in \liealg g_{\field C}^\ast$, we still have that }\begin{equation*}\tag{$\ast$}
        \eta\wedge (\DECC \eta)^{k+1}=0.
    \end{equation*}
    Since extending forms $\mathbb{C}$-multi-linearly yields a map of differential graded commutative algebras
    \[(\wedge^{\bullet} \mathfrak{g}^*,\wedge,\DEC)\to (\wedge_{\mathbb{C}}^{\bullet} \mathfrak{g}_{\mathbb{C}}^*,\wedge,\DECC),\]
    equation $(\ast)$ holds for elements of $\liealg g^\ast$. Thus, it is enough to show that, if $\nu,\mu\in \liealg g^\ast$, then also $\nu+ i\mu$ satisfies $(\ast)$. 
    Consider the map
    \begin{equation*}
        \begin{aligned}
            q\colon \field C&\to \wedge_{\mathbb{C}}^{2k+3}\liealg g_{\field C}^\ast \\
            z&\mapsto (\nu+z\mu)\wedge (\DECC (\nu+z\mu))^{k+1}.
        \end{aligned}
    \end{equation*}
    Then $q$ is a complex polynomial that vanishes for all $z\in \field R$, hence $q$ vanishes identically, proving the claim.

    Lemma \ref{lemma:equivalence_orbit_dga} clearly also holds for complex Lie algebras. By $(\ast)$, any element $\eta\in U$ is of type \ref{type_2}. Since this is equivalent to \emph{(2'')}, there exists $Y\in \mathfrak{g}_{\mathbb{C}}$ such that
\begin{equation}\label{eq:radial_tangent}
        \eta=\mathrm{ad}^*_Y\eta.
    \end{equation}    
Let $B\colon \liealg g_{\field C}\times\liealg g_{\field C} \to \field C$ denote the Killing form of $\liealg g_{\field C}$. Then, for $\eta=B^{\sharp}X$, \eqref{eq:radial_tangent} is equivalent to 
    \[X=[X,Y].\]
Indeed, this follows because  $B^{\sharp}$ is equivariant, and so 
\[\mathrm{ad}_Y^*B^{\sharp}X=-B^{\sharp}(\mathrm{ad}_Y(X))=B^{\sharp}([X,Y]).\]

    On the other hand, by \cite[\S II, 2]{knapp:2002a} the set of regular elements in $\liealg g_{\field C}$ is given by the non-vanishing of a polynomial, hence it is also open and dense. Therefore, there exists a regular $X\in \liealg g_{\field C}$ such that $B^\sharp X\in U$. Since $X$ is regular, by \cite[\S II, Theorem 2.9]{knapp:2002a},
    \begin{equation*}
        \liealg h:=\{ H\in \liealg g_{\field C}\colon [X, H]=0 \}
    \end{equation*}
is a Cartan subalgebra. Consider the corresponding root space decomposition
    \begin{equation*}
        \liealg g_{\field C}=\bigoplus_{\alpha\in \Phi} \liealg g_\alpha\oplus \liealg h.
    \end{equation*}
Since $B^{\sharp}X\in U$, there exists 
    \begin{equation*}
        Y=\sum_{\alpha\in \Phi}Y_\alpha+Y_0\in \mathfrak{g}_{\mathbb{C}},
    \end{equation*}
    such that $X=[X,Y]$.
Since $X\in \mathfrak{h}$, we obtain a contradiction: 
    \begin{equation*}
        X=[X,Y]=\sum_{\alpha\in \Phi}\alpha(X)Y_\alpha\in \mathfrak{h}\cap \big(\bigoplus_{\alpha\in \Phi} \liealg g_\alpha\big) =\{0\},
    \end{equation*}
which concludes the proof.
\end{proof}

\begin{lemma}\label{lemma:ss_constant_height_compact}
    A semisimple Lie algebra of constant height is compact.
\end{lemma}
\begin{proof}
    Let $\liealg g$ be a semisimple Lie algebra of constant height. Let $\theta\colon \liealg g\to \liealg g$ be a Cartan involution with corresponding Cartan decomposition $\liealg g=\liealg k\oplus \liealg p$, where
    \begin{equation*}
            \liealg k=\{ X\in \liealg g\colon \theta(X)=X \}, \quad\qquad
            \liealg p=\{ X\in \liealg g\colon \theta(X)=-X \}.
    \end{equation*}
 If $B$ denotes the Killing form, then the two-form 
    \begin{equation*}
        \begin{aligned}
            B_\theta(X,Y):=-B(X,\theta(Y)), \qquad X,Y\in \mathfrak{g},
        \end{aligned}
    \end{equation*}
    is positive definite and symmetric, and for any $Y\in \liealg p$ the operator $\ad_Y\colon \liealg g\to \liealg g$ is self-adjoint, see \cite[Lemma 6.27]{knapp:2002a}. 
    Thus, $\ad_Y$ is diagonalisable over $\field R$. Let $X\in \liealg g$ be an eigenvector of $\ad_Y$ with eigenvalue $\lambda\in \field R\setminus \{0\}$. 
    Then
    \begin{equation*}
        [X,-\tfrac{1}{\lambda}Y]=X.
    \end{equation*}
    As we have seen in the proof of Lemma \ref{lemma:ss_constant_height_is_type_1}, this is equivalent to $B^{\sharp}(X)$ being of type \ref{type_2}. This contradicts Lemma \ref{lemma:ss_constant_height_is_type_1}.
    Hence, all eigenvalues of $\ad_Y$ are zero. Since $\ad_Y$ is diagonalizable, $Y$ must lie in the centre of $\liealg g$. 
    Since $\liealg g$ is semisimple, $Y=0$. We have shown that $\liealg p=\{0\}$, which by \cite[Theorem 3.6.2]{duistermaat.kolk:2000a} implies that $\liealg g$ is compact, as its Killing form is negative definite. 
\end{proof}

\begin{theorem}\label{theorem:ss_constant_height_is_so3}
Any semisimple Lie algebra of constant height is isomorphic to $\liealg{so}(3)$.
\end{theorem}
\begin{proof}
Let $\mathfrak{g}$ be semisimple of constant height $k$. By Lemma \ref{lemma:ss_constant_height_is_type_1}, all elements of $\liealg g^\ast$ have type \ref{type_1}. Hence, by Lemma \ref{lemma:equivalence_orbit_dga}, all non-trivial coadjoint orbits of $\liealg g^\ast$ have dimension $2k$. The Killing form induces an isomorphism $\liealg g^\ast \simeq \liealg g$ of $\liealg g$-representations. So in the notation of \cite[Definition (2.8.3)]{duistermaat.kolk:2000a},
    \begin{equation*}\tag{$\ast$}
        \liealg g^\mathrm{reg}=\liealg g\setminus \{0\}.
    \end{equation*}
    Let $\liealg t\subseteq \liealg g$ be a maximal abelian subalgebra. 
    Since by Lemma \ref{lemma:ss_constant_height_compact} $\liealg g$ is compact, the complexification $\liealg t_{\field C}$ is a Cartan subalgebra of $\liealg g_{\field C}$, see \cite[Proposition 6.47]{knapp:2002a}, and corresponding roots $\alpha\in \Phi$ take real values on $i\liealg t$, see \cite[Corollary 4.49]{knapp:2002a}.
    By \cite[Theorem (3.7.1) (ii)]{duistermaat.kolk:2000a} we have
    \begin{equation*}
        \liealg t\setminus \{0\} \stackrel{(\ast)}{=}\liealg g^\mathrm{reg}\cap \liealg t=\liealg t\setminus\bigcup_{\alpha\in \Phi}\ker\alpha.
    \end{equation*}
    Since $\alpha\colon i\liealg t\to \field R$, this implies $\dim \liealg t=1$, i.e.\ $\liealg g_{\field C}$ has rank 1 and is thus isomorphic to $\liealg{sl}_2(\field C)$. Its compact real form is $\liealg{so}(3)$, which proves the statement.    
\end{proof}

With Theorem \ref{theorem:ss_constant_height_is_so3} at hand, we can proceed with the proof of Theorem \ref{theorem:classification_LA_of_constant_height}. 
First, we show that constant height Lie algebras have to be abelian extensions of $\field R$ or $\liealg{so}(3)$.

\begin{lemma}\label{lemma:ideals}
    Let $\liealg g$ be a Lie algebra of constant height $k$ and $\liealg h\subsetneq \liealg g$ a proper ideal. Then: 
    \begin{enumerate}
        \item The height of $\liealg g /\liealg h$ is constant, equal to that of $\liealg g$.
        \item The ideal $\liealg h$ is abelian.
    \end{enumerate}
\end{lemma}
\begin{proof}
    The first part follows from the fact that the identification 
    \begin{equation*}
        (\liealg g /\liealg h)^\ast=\liealg h^\circ \subseteq \liealg g^\ast
    \end{equation*}
    induces an inclusion of differential graded commutative algebras 
    \begin{equation*}
        \big(\wedge^\bullet (\liealg g /\liealg h)^\ast,\wedge , \mathrm{d}_{\mathfrak{g}/\mathfrak{h}}) \hookrightarrow \big(\wedge^\bullet \liealg g^\ast, \wedge, \DEC\big).
    \end{equation*}
    For the second part, choose a complement $\liealg g=\liealg h\oplus\liealg c$. For all $\theta\in \liealg h^\circ\setminus\{0\}$ and $\xi\in \liealg h^\ast=\liealg c^\circ$ we have
    \begin{equation*}
        (\theta+\xi)\wedge(\DEC(\theta+\xi))^{k+1}=0.
    \end{equation*}
    For a fixed $\theta$, this equation is a polynomial in (the components of) $\xi$.
    Thus, to vanish, each homogeneous component needs to vanish separately. 
    The linear terms read
    \begin{equation*}
        \xi\wedge (\DEC \theta)^{k+1}+(k+1)\theta\wedge(\DEC \theta)^k\wedge \DEC \xi.
    \end{equation*}
    Using that $\DEC \theta\in \wedge^2 \mathfrak{h}^{\circ}$, applying $\mathrm{i}_X \mathrm{i}_Y $ to this equation, where $X,Y\in \liealg h$, gives
    \begin{equation*}
        (k+1)\xi([X,Y])\theta\wedge(\DEC \theta)^k=0.
    \end{equation*}
    Since this holds for all $\xi\in \liealg h^\ast$, the assumption $\theta\wedge(\DEC \theta)^k\neq 0$ implies $[X,Y]=0$.
\end{proof}

\begin{corollary}\label{corollary:constant_height_abelian_ext_of_so3}
   A Lie algebra of constant height has height either $0$ or $1$.
\end{corollary}
\begin{proof}
    Let $\liealg h\subsetneq \liealg g$ be a maximal proper ideal (which is necessarily abelian by Lemma \ref{lemma:ideals}). Then $\liealg g /\liealg h$ is a Lie algebra of constant height with no proper ideals, so either $\liealg g /\liealg h$ is simple or $\liealg g /\liealg h=\field R$. In the latter case the height of $\liealg g$ clearly is $0$. If $\liealg g /\liealg h$ is simple, then $\liealg g /\liealg h\simeq \liealg{so}(3)$ by Theorem \ref{theorem:ss_constant_height_is_so3}.
\end{proof}

\begin{proposition}
A Lie algebra of constant height $0$ is either abelian or is isomorphic to $\field R\ltimes \mathbb{R}^n$, for the diagonal representation of $\field R$ on $\mathbb{R}^n$. 
\end{proposition}
\begin{proof}
Let $\mathfrak{h}\subsetneq \mathfrak{g}$ be a maximal proper ideal. By Lemma \ref{lemma:ideals}, $\mathfrak{h}$ is abelian and $\liealg g/\liealg h$ is one-dimensional. 
 Let $e\in \liealg g/ \liealg h\setminus \{0\}$ and $X\in \liealg g$ be a preimage of $e$ under the quotient map. Then we can identify $\mathfrak{g}= \mathbb{R}X\oplus\mathfrak{h}$, and the Lie bracket on $\liealg g$ is completely determined by the linear map
    \begin{equation*}
        A:=[X,\argument]\colon \liealg h\to \liealg h.
    \end{equation*}
   Let $\theta\in (\field R X)^\circ=\liealg h^\ast$ be given. By the assumption on $\liealg g$, we have $\theta\wedge \DEC \theta=0$. Inserting $X$ yields
    \begin{equation*}
        0=\mathrm{i}_X (\theta \wedge \DEC \theta)=-\theta\wedge\mathrm{i}_X\DEC \theta = \theta\wedge  A^\ast \theta.
    \end{equation*}
    Hence, $A^\ast\theta=\lambda_\theta \theta$, for some $\lambda_{\theta}\in \mathbb{R}$. Therefore, any vector of $\liealg h^\ast$ is an eigenvector of $A^\ast$. Thus, the eigenvalues have to coincide, i.e.\ there exists $\lambda\in \mathbb{R}$, such that 
    \begin{equation*}
        A^\ast = \lambda \id_{\liealg h^\ast} \Leftrightarrow A=\lambda\id_{\liealg h}.
    \end{equation*}
    If $\lambda=0$, then $\mathfrak{g}$ is abelian. If $\lambda \neq 0$, by rescaling $X$, we obtain the isomorphism $\liealg g\simeq \field R\ltimes \liealg h$, for the diagonal representation of $\field R$ on $\liealg h$.
\end{proof}

Finally, we show that non-trivial abelian extensions of $\liealg{so}(3)$ do not have constant height.

\begin{proposition}
Any Lie algebra of constant height $1$ is isomorphic to $\liealg{so}(3)$. 
\end{proposition}
\begin{proof}
 Let $\liealg g$ be a Lie algebra of constant height $1$.
    By the proof of Corollary \ref{corollary:constant_height_abelian_ext_of_so3} there is a short exact sequence
    \begin{equation}\label{eq:sesq_1}
        0\to\liealg h\to\liealg g\to \liealg{so}(3)\to 0,
    \end{equation}
    where $\liealg h$ is abelian. 
    The adjoint representation of $\mathfrak{g}$ restricted to $\mathfrak{h}$ descends to a representation $\rho\colon \liealg{so}(3)\to \liealg{gl}(\liealg h)$. 
    
    The obstruction to find a splitting $\sigma\colon\mathfrak{so}(3)\to \mathfrak{g}$ of \eqref{eq:sesq_1} compatible with Lie brackets lies in the second cohomology group of $\liealg{so}(3)$ with coefficients in $\mathfrak{h}$. By the Whitehead Lemma, $\mathrm{H}^2(\liealg{so}(3), \mathfrak{h})=0$. Thus, we may assume that $\liealg g$ is a semidirect product   
    \begin{equation*}
        \liealg g=\liealg{so}(3)\ltimes \liealg h.
    \end{equation*}
    We show that the representation $\rho$ of $\liealg{so}(3)$ on $\mathfrak{h}$ is trivial.

    {\sf Claim:} {\it For all $\xi\in \liealg h^\ast$ we have
    \begin{equation*}
        \dim \Span \{ \xi, \rho_X^\ast\xi\, \colon\,  X\in\liealg{so}(3)\} \leq 2.  
    \end{equation*}
    }
    Fix $\theta\in \liealg{so}(3)^\ast$. As in the proof of Lemma \ref{lemma:ideals}, the expression
    \begin{equation*}
        (\theta+\xi)\wedge(\DEC(\theta+\xi))^{2}=0
    \end{equation*}
    is a polynomial in $\xi\in \liealg h^\ast$. 
    Vanishing of the degree two component reads
    \begin{equation*}\tag{$\ast$}
        \theta\wedge (\DEC\xi)^2+2\xi\wedge \DEC \theta\wedge \DEC\xi=0.
    \end{equation*}
    We evaluate this expression on elements of $\liealg{so}(3)$. Denote by $\{X_1,X_2,X_3\}$ the standard basis of $\liealg{so}(3)$ and by $\{\theta_1,\theta_2,\theta_3\}$ the dual basis. Then
    \begin{equation*}
        \begin{aligned}
            \mathrm{i}_{X_i}\DEC \theta_j=-\theta_j([X_i, \argument])=\sum_{k=1}^3\varepsilon_{ijk}\theta_k
        \end{aligned}
    \end{equation*}
    and, for $\xi\in \liealg h^\ast$, 
    \begin{equation*}
        \mathrm{i}_X\DEC \xi=-\xi([X,\argument])=-\xi(\rho_X(\argument))=-\rho_X^\ast\xi.
    \end{equation*}
    We take $(\ast)$ for $\theta=\theta_1$ and evaluate on $X_1$, $X_2$, and $X_3$:
    \begin{equation*}
        \begin{aligned}
            0&=\mathrm{i}_{X_3}\mathrm{i}_{X_2}\mathrm{i}_{X_1}\big(\theta_1\wedge (\DEC\xi)^2+2\xi\wedge \DEC \theta_1\wedge \DEC\xi\big)\\
            &=\mathrm{i}_{X_3}\mathrm{i}_{X_2}\big( (\DEC\xi)^2+2 \theta_1\wedge \rho_{X_1}^\ast\xi\wedge \DEC\xi+2\xi\wedge \DEC \theta_1\wedge \rho_{X_1}^\ast\xi\big)\\
            &=\mathrm{i}_{X_3} ( -2 \rho_{X_2}^\ast\xi \wedge \DEC \xi +2\theta_1\wedge\rho_{X_1}^\ast\xi\wedge \rho_{X_2}^\ast\xi +2\xi\wedge\theta_3\wedge\rho_{X_1}^\ast\xi ) \\
            &=- 2\rho_{X_2}^\ast\xi \wedge \rho_{X_3}^\ast\xi-2\xi\wedge \rho_{X_1}^\ast\xi.
        \end{aligned}
    \end{equation*}
    In conclusion, for all $\xi\in \liealg h^\ast$, 
    \begin{equation*}
        \xi\wedge\rho_{X_1}^\ast\xi=\rho_{X_3}^\ast\xi \wedge \rho_{X_2}^\ast \xi.
    \end{equation*}
    Hence, for any three elements $\eta_1,\eta_2,\eta_3$ of the set 
    \begin{equation*}
        \{ \xi, \rho_{X_1}^\ast\xi ,\rho_{X_2}^\ast\xi , \rho_{X_3}^\ast \xi \}    
    \end{equation*}
    we have $\eta_1\wedge\eta_2\wedge\eta_3=0$, proving the claim. 

    With this, we show that the representation $\rho$ is actually trivial. Consider the integration of $\rho$ to an action of the compact group $SU(2)$ on $\liealg h^\ast$. By compactness of $SU(2)$ there exists an invariant inner product on $\liealg h^\ast$ with induced norm $\norm{\argument}$. Let $\xi\in \liealg h^\ast$ and $r:=\norm{\xi}$. Let $\mathcal{O}_\xi$ denote the orbit of $\xi$. Then, since $\mathcal{O}_\xi\subseteq S_r(0)$ is contained in a sphere,
    \begin{equation*}
        \field R \xi \cap T_\xi \mathcal{O}_\xi=\{0\}.
    \end{equation*}
    Using 
    \begin{equation*}
        T_\xi \mathcal{O}_\xi=\{ \rho_{X}^\ast\xi \colon X\in \liealg{so}(3)\},
    \end{equation*}
    the Claim implies that $\dim T_\xi \mathcal{O}_\xi\leq 1$. Suppose that $\dim  T_\xi \mathcal{O}_\xi=1$. 
    Then the isotropy Lie algebra of $\xi$ is a two dimensional subalgebra of $\liealg{so}(3)$. But there exists no such subalgebra.
    Hence, $\dim T_\xi \mathcal{O}_\xi=0$ and the representation has to be trivial.

    In conclusion, $\liealg g=\liealg{so}(3)\oplus \liealg h$ is a direct sum of Lie algebras. If $\liealg h\neq 0$, then $\liealg{so}(3)\subsetneq \liealg g$ is a proper ideal. By Lemma \ref{lemma:ideals} $\liealg{so}(3)$ is abelian, a contradiction. Hence, $\liealg h=0$ and $\liealg g\simeq \liealg{so}(3)$. 
\end{proof}

\bibliographystyle{alpha}


\end{document}